\documentclass[oneside,11pt]{amsart}

\usepackage{float}
\usepackage{amscd,amssymb}
\usepackage[mathscr]{eucal}
\usepackage[all]{xy}
\usepackage{mathrsfs}
\usepackage{amsfonts}
\usepackage{amsmath}
\usepackage{amsthm}
\usepackage{latexsym}
\usepackage[all]{xy}

\title{Tilting bundles via the Frobenius morphism}

\author{Alexander Samokhin}

\address{Institute for Information Transmission Problems, Moscow,
  Russia}

\email{alexander.samokhin@gmail.com}

\newcommand{\Oo}{\mathcal O}
\newcommand{\Uu}{\mathcal U}

\newcommand{\Pp}{\mathbb P}

\newcommand{\Ff}{\mathcal F}
\newcommand{\Ee}{\mathcal E}
\newcommand{\Ll}{\mathcal L}

\newcommand{\D}{\mathcal D}

\newcommand{\T}{\mathcal T}

\newcommand{\V}{\sf V}
\newcommand*{\RHom}{\mathop{\mathrm RHom}\nolimits}
\newcommand*{\Hom}{\mathop{\mathrm Hom}\nolimits}
\newcommand*{\Dd}{\mathop{\mathrm D\kern0pt}\nolimits}
\newcommand*{\DD}{\mathop{\mathbb D\kern0pt}\nolimits}

\newcommand*{\Ext}{\mathop{\mathrm Ext}\nolimits}
\newcommand{\fmod}{{\amod^{\text{{\tt\tiny fg}}}}}
\newcommand{\amod}{{\text{\rm -mod}}}

\newtheorem{theorem}{Theorem}[section]
\newtheorem{corollary}{Corollary}[section]
\newtheorem{lemma}{Lemma}[section]
\newtheorem{proposition}{Proposition}[section]

\newtheorem{remark}{Remark}[section]
\newtheorem{definition}{Definition}[section]

\numberwithin{equation}{section}

\long\def\comment#1{}

\begin{document}


\maketitle






\begin{abstract}
Let $X$ be a smooth algebraic variety over an algebraically closed
field $k$ of characteristic $p>0$, and ${\sf F}\colon X\rightarrow X$ 
the absolute Frobenius morphism. In this paper we compute the cohomology groups ${\rm H}^{i}(X,{\mathcal End}({\sf
  F}_{\ast}\Oo _X))$ and show that these groups vanish for $i>0$ in a
number of cases that include some toric Fano varieties, blowups of the
projective plane (e.g., Del Pezzo surfaces), and the homogeneous
spaces ${\bf G}/{\bf P}$ of groups of type ${\bf A}_2$ and ${\bf
  B}_2$. The reason to consider such cohomology
groups is that if the higher cohomology vanishing holds, then the sheaf ${\sf F}_{\ast}\Oo _X$
gives a tilting bundle on $X$, provided that it is a
generator in $\Dd ^{b}(X)$, the bounded derived category of coherent
sheaves on $X$. For homogeneous spaces ${\bf G}/{\bf P}$ the derived localization theorem \cite{BMR} implies
that the sheaf ${\sf F}_{\ast}\Oo _{{\bf G}/{\bf P}}$ is a generator
if $p>h$, the Coxeter number of the group $\bf G$, thus producing a
tilting bundle on ${\bf G}/{\bf P}$ for the group $\bf G$ of type
${\bf A}_2$ and ${\bf B}_2$. We also verify that ${\sf F}_{\ast}\Oo _X$ is a generator for several
non-homogeneous cases.
\end{abstract}


\section*{Introduction}


Let $X$ be a smooth algebraic variety and $\Dd ^{b}(X)$ its bounded
derived category of coherent sheaves. The now classical theorem of
Beilinson \cite{Beil} states that if $X= \Pp ^n$, then there is an equivalence of categories:
\begin{equation}\label{eq:tiltequiv} 
\Phi \colon \Dd ^{b}(X) \simeq \Dd ^{b}({\sf A}-\mbox{mod})\nonumber ,
\end{equation}

\noindent where $\Dd ^{b}({\sf A}-\mbox{mod})$ is the bounded
derived category of finitely generated left modules over a 
non-commutative algebra $\rm A$, the path algebra of Beilinson's
quiver. The proof consists of finding a vector bundle $\Ee$ on $\Pp ^n$ such
that the bundle of endomorphisms of $\Ee$ does not have higher
cohomology groups and the direct summands of $\Ee$ generate the whole
$\Dd ^{b}(\Pp ^n)$. Such a bundle is called a tilting generator, and 
equivalences of the above type are called tilting equivalences. The
direct summands of $\Ee$ form a so-called full strong exceptional collection
\cite{Bo}. Therefore, given a variety $X$, one way to construct a tilting equivalence is to
find such a collection on $X$.

Some time later tilting equivalences appeared in the proof of
the derived McKay correspondence \cite{BKR} and were subsequently interpreted,
for a broader class of examples, as non-commutative resolutions of
singularities \cite{VdB}. In particular, tilting equivalences 
were shown to exist for small contractions of 3-folds ({\it loc.cit.}) and, more
recently, for more general contractions \cite{TU}. 

Recent advances in representation theory in positive characteristic
\cite{BezICM} provide us with another rich source of tilting
equivalences for varieties equipped with a non-degenerate algebraic
symplectic 2-form. On one hand, these examples fit into the framework of non-commutative resolution of
singularities. On the other hand, in the core of these equivalences
are specific features of positive characteristic. In a nutshell, tilting bundles obtained this way arise 
as splitting bundles of some Azumaya algebras. The latter turn out to be non-commutative deformations
(quantizations) of the sheaf of functions on the Frobenius twist of a symplectic variety in
positive characteristic. This was used, for instance, in \cite{SymMcKay} to prove 
the derived McKay correspondence in the symplectic case. The general
framework of quantization of symplectic varieties in positive
characteristic was laid down in \cite{BezKal}. Applications of this
theory to tilting equivalences for general symplectic resolutions were
given in \cite{Kaledquant}.

Cotangent bundles of flag varieties ${\bf G}/{\bf B}$ are a prominent example of this theory. Their quantization in positive characteristic leads to the
derived Beilinson--Bernstein localization theorem -- the subject of the
influential work \cite{BMR}. The non-commutative sheaf of algebras that one
obtains quantizing the cotangent bundle of an algebraic variety $X$ is
the sheaf of $\rm PD$-differential operators, whose central reduction
is the sheaf of small differential operators ${\mathcal End}({\sf F}_{\ast}\Oo
_X)$ (a split Azumaya algebra), ${\sf F}$ being the Frobenius morphism. An immediate consequence of the derived
localization theorem \cite{BMR} is that, for sufficiently large $p$, the bundle ${\sf F}_{\ast}\Oo
_{{\bf G}/{\bf B}}$ generates the derived category of coherent sheaves
on  ${\bf G}/{\bf B}$. More generally, under similar assumptions on
$p$, the bundle ${\sf F}_{\ast}\Oo
_{{\bf G}/{\bf P}}$ generates the derived category of ${\bf G}/{\bf
  P}$ for a parabolic subgroup ${\bf P}\subset {\bf G}$ (see Corollary
\ref{cor:Frobgeneratorforp>h}). 

Inspired by these results, we would like to apply the characteristic
$p$ methods to find tilting bundles. 
In the present paper, however, we are mostly interested in the proper and non-symplectic case,
while the quantization theory in positive characteristic works well
for non-compact symplectic varieties \cite{Kaledquant}. In the realm
of smooth proper varieties, the most classical examples from the point of view of tilting equivalences
are homogeneous spaces ${\bf G}/{\bf P}$ and toric varieties. Strong exceptional collections (hence,
tilting bundles) on some homogeneous spaces have been known since the seminal work of Kapranov
\cite{Kap}. It has since been conjectured by various authors that homogeneous spaces of
any semisimple algebraic group should have full -- possibly, strong --
exceptional collections (see \cite{Kuz} for a recent progress via the
theory of exceptional collections). Similarly, King, in the
well-known unpublished preprint \cite{King}, has conjectured that a
smooth complete toric variety should have a
tilting bundle, whose summands would have been line bundles. In such
a generality this conjecture has turned out to be false, as was recently shown by Hille and Perling
\cite{Kingcountex}. On the other hand, many examples of
toric varieties that have strong exceptional collections of line
bundles \cite{HP} indicate that a refinement of King's conjecture
could still be true.

Thus, given a variety $X$ that conjecturally has a tilting bundle,
let us reduce $X$ to positive characteristic and consider the sheaf
${\sf F}_{\ast}\Oo _X$ (more generally, the sheaf ${\sf F}^n_{\ast}\Ll$
for a line bundle $\Ll$ on $X$ for iterated Frobenius morphism). If
$X={\bf G}/{\bf P}$ is a homogeneous space and $p$ is large enough, then, by Corollary
\ref{cor:Frobgeneratorforp>h}, the vanishing of higher cohomology
groups ${\rm H}^{i}({\bf G}/{\bf P},{\mathcal End}({\sf F}_{\ast}\Oo
_{{\bf G}/{\bf P}}))$ implies that the bundle ${\sf F}_{\ast}\Oo _{{\bf G}/{\bf
    P}}$ is tilting on ${\bf G}/{\bf P}$.  Independently, such a vanishing appeared earlier in a different context,
namely in the study of ${\sf D}$-affinity of flag varieties in
positive characteristic \cite{Haas}, \cite{AK} (in the next paragraph we elaborate on
these results). Hence, there is evidence showing that it may be sensible to
compute the Frobenius pushforward of the structure sheaf on ${\bf G}/{\bf P}$. This was the
motivation to also consider the Frobenius pushforwards of line bundles
on non-homogeneous varieties in hope to get tilting bundles -- for
example, for toric varieties from the previous paragraph.
Considering the Frobenius pushforwards of line bundles on toric
varieties seems to be relevant to our problem thanks to Thomsen's theorem \cite{Bogv} that asserts that the Frobenius pushforward of a line bundle on a toric variety is the
direct sum of line bundles. Though it is not true in general that the
Frobenius pushforward of the structure sheaf of a toric variety is tilting, quite a few examples of
tilting bundles can be obtained along this way (see Section
\ref{sec:toricFanovar}). Note that, in characteristic zero,
multiplication maps on toric varieties were used in \cite{Botoric} to
relate the derived category of coherent sheaves and the derived
category of constructible sheaves on corresponding real tori, as
predicted by the homological mirror symmetry. The desired equivalence
of categories was established in \cite{BCTZ}. From the
${\sf D}$--affinity perspective, it is known, however, that the only $\sf
D$-affine smooth toric varieties are products of projective spaces \cite{ThD-afftoric}. 

The paper is organized as follows. We begin with introductory Section
\ref{sec:prelim} where we recall the necessary facts about the Frobenius morphism, differential
operators, vanishing theorems for line bundles on flag varieties, and tilting bundles and
derived equivalences. We formulate the derived localization theorem \cite{BMR} that implies 
that the bundle ${\sf F}_{\ast}\Oo _{{\bf G}/{\bf P}}$ is a generator
in the derived category (Corollary \ref{cor:Frobgeneratorforp>h}). 
Section \ref{subsec:genstat} contains a few technical statements
that are needed to compute the groups $\Ext ^{i}({\sf F}_{\ast}\Oo
_X,{\sf F}_{\ast}\Oo _X)$. In particular, we prove short exact
sequences that help to compute these groups when the variety in
question is either a $\Pp ^1$-bundle over a smooth base or the blow-up
of a smooth surface at a number of points in general position. In Section
\ref{sec:flagvar} we give first applications of the techniques
developed in the previous sections. We prove in Theorems
\ref{th:SL_3} and \ref{th:B_2} that ${\rm H}^{i}({\bf G}/{\bf B}, {\mathcal End}({\sf F}_{\ast}\Oo _{{\bf
    G}/{\bf B}})) = 0$, if $\bf G$ is of type ${\bf A}_2$ or ${\bf
  B}_2$ (in the latter case $p$ is odd). Note that the above vanishing
has previously been known in both cases: in type ${\bf A}_2$
this was proved by Haastert in \cite{Haas} and in type ${\bf B}_2$
by Andersen and Kaneda in \cite{AK} (for the case $p=2$ as well,
while we have to restrict ourselves to odd $p$ to be able to use the
Kumar--Lauritzen--Thomsen vanishing theorem; however, to obtain a
tilting bundle out of ${\sf F}_{\ast}\Oo _{{\bf G}/{\bf B}}$ the
characteristic of $k$ must anyway be greater than 3). The method used in \cite{Haas} and \cite{AK} was to identify the sheaf  
${\mathcal End}({\sf F}_{\ast}\Oo _{{\bf G}/{\bf B}})$ with an
equivariant vector bundle on ${\bf G}/{\bf B}$ associated to the
induced module ${\rm Ind}_{\bf B}^{{\bf G}_1{\bf B}}(2(p-1)\rho)$
(called the Humphreys--Verma module in \cite{KanHVmod}), where ${\bf G}_1$ is the first Frobenius kernel, and to study an
appropriate filtration on such a module to find the
weights of composition factors. Our proof is different and based on the properties of the algebra of crystalline
differential operators, one of the main ingredients in \cite{BMR},
that help to reduce the problem to the cohomology groups of the Frobenius pullback
of a certain Koszul complex. These groups can be computed using vanishing theorems for line bundles on ${\bf G}/{\bf B}$ (and on
the total spaces of cotangent bundles of these) and techniques
developed in Section \ref{subsec:genstat}. In a subsequent version of
the paper we show, using the above methods, that for sufficiently large $p$ the Frobenius
pushforward of the structure sheaf is tilting on the flag variety
${\bf G}/{\bf B}$, where the root system of $\bf G$ is of type ${\bf
  G}_2$.

In a companion paper \cite{Sam2} we apply similar arguments to quadrics
  and some partial flag varieties of type ${\bf
  A}_n$ to prove that the Frobenius pushforward of the structure sheaf
  gives a tilting bundle on these varieties. Further examples will appear in \cite{Sam3}. 
We hope that this approach can help to compute the higher Frobenius
  pushforwards of structure sheaves and to give an insight into the $\sf D$-affinity of flag varieties in positive characteristic. For
instance, the argument used in the proof of Theorem \ref{th:SL_3} can easily be
extended to show the $\sf D$-affinity of the flag variety in type
${\bf B}_2$ (see \cite{SamD-affflag}).

In Section \ref{sec:toricFanovar} we work out several examples of
toric Fano varieties. In particular, we study the vanishing behaviour
of the groups $\Ext ^{i}({\sf F}_{\ast}\Oo
_X,{\sf F}_{\ast}\Oo _X)$ for twelve out of a total of eighteen smooth toric Fano
3-folds \cite{Ba}; it turns out that these groups vanish for $i>0$ except
for two cases. On the other hand, it is easy to check, in the case
when the vanishing holds, that the bundle ${\sf F}_{\ast}\Oo _X$ is a generator,
hence tilting. Using \cite{Bogv}, it should be possible 
to compute these groups for all smooth toric Fano 3-folds. Finally, in Section
\ref{sec:blowups} we study the bundle ${\sf F}_{\ast}\Oo _{X_k}$,
where $X_k$ is the blow-up of $\Pp ^2$ at $k$ points in general
position, and show that $\Ext ^{i}({\sf F}_{\ast}\Oo _{X_k},{\sf
  F}_{\ast}\Oo _{X_k}) = 0$ for $i>0$. 

\subsection*{Acknowledgements}
The first draft of this paper with most of the results presented here was
written in the spring of 2005. The present work was done during the author's stays at 
the MPI Bonn in the spring of 2007 (see \cite{SamMPI}), and at the IHES and the Institute
of Mathematics of Jussieu in the spring of 2009. He gratefully
acknowledges their hospitality and support.

\subsection*{Notation} Throughout the paper the ground field is
an algebraically closed field $k$ of characteristic
$p>0$. A variety means an integral separated scheme of finite
type over $k$. For a variety $X$ its bounded derived
category of coherent sheaves is denoted ${\rm D}^{b}(X)$, and $[1]$
denotes the shift functor. For an object $\Ee$ of $\Dd ^{b}(X)$ denote
$\Ee ^{\vee}$ the dual to $\Ee$, that is the object ${\mathcal
  RHom}^{\bullet}(\Ee ,\Oo _X)$. The symbol $\boxtimes$ denotes
external tensor product.

\tableofcontents


\section{Preliminaries}\label{sec:prelim}


\subsection{The Frobenius morphism (\cite{Il})}\label{subsec:Frobmorph}

Let $X$ be a smooth variety over $k$. The absolute Frobenius morphism ${\sf F}_X$ is an endomorphism
of $X$ that acts identically on the topological space of $X$ and
raises functions on $X$ to the $p$-th power:
\begin{equation}\label{eq:defofFrobmor}
{\sf F}_X\colon X\rightarrow X,\qquad f\in \Oo _X\rightarrow f^p\in \Oo _X.
\end{equation}

Let $\pi\colon X\rightarrow S$ be a morphism of
$k$-schemes. Then there is a commutative square:

\begin{figure}[H]
$$
\xymatrix @C5pc @R4pc {
 X\ar[r]^{{\sf F}_X} \ar@<-0.1ex>[d]^{\pi} & X \ar@<-0.4ex>[d]^{\pi} \\
       S   \ar[r]^{{\sf F}_S} & S
 }
$$
\end{figure}

Denote $X'$ the scheme $(S,{\sf F}_S)\times _S X$ obtained
by the base change under ${\sf F}_S$ from $X$. The morphism ${\sf
  F}_X$ defines a unique $S$-morphism ${\sf F}={\sf F}_{X/S}\colon
  X\rightarrow X'$, such that there is a commutative diagram:

\begin{figure}[h]\label{fig:fig14}
$$
\xymatrix @C5pc @R4pc {
 X\ar[r]^{\sf F_{X/S}} \ar[dr]^{\pi} & X'\ar[r]^{\phi}
 \ar[d]_{\pi '} & X \ar[d]_{\pi}\\
& S\ar[r]^{\sf F_{S}} & S
 }
$$
\end{figure}

The composition of upper arrows $\phi \circ {\sf F}$ is
equal to ${\sf F}_X$, and the square is cartesian. The morphism ${\sf
  F}$ is said to be the relative Frobenius morphism of $X$ over
  $S$. The morphism ${\sf F}_X$ is not a morphism of $S$-schemes. On
  the contrary, the morphism ${\sf F}_{X/S}$ is a morphism of $S$-schemes.

\begin{proposition}\label{prop:Frobprop}
Let $S$ be a scheme over $k$, and $\pi \colon X\rightarrow S$ a smooth
morphism of relative dimension $n$. Then the relative Frobenius
morphism ${\sf F}\colon X\rightarrow X'$ is a finite flat morphism,
and the $\Oo _X'$-algebra ${\sf F}_{\ast}\Oo _X$ is locally free of
rank $p^n$.
\end{proposition}

Let $\pi \colon X\rightarrow S$ be a smooth morphism as
above, and $\Ff$ a coherent sheaf (a complex of coherent sheaves) on $X$. 
Proposition \ref{prop:Frobprop} implies:
\begin{equation}\label{eq:Frobcommutativitywithmorphisms}
{\rm R}^{i}\pi _{\ast}{{\sf F}_{X}}_{\ast}\Ff = {{\sf F}_S}_{\ast}{\rm
  R}^{i}\pi _{\ast}\Ff.
\end{equation}

Indeed, it follows from the spectral sequence for the composition of two
  functors:
\begin{equation}
{\rm R}^{i}\pi _{\ast}{{\sf F}_{X}}_{\ast}\Ff = {\rm R}^{i}(\pi \circ
{\sf F}_X)_{\ast}\Ff = {\rm R}^i({\sf F}_S \circ \pi)_{\ast}\Ff =
{{\sf F}_S}_{\ast}{\rm R}^{i}\pi _{\ast}\Ff.
\end{equation}

Let $S = {\rm Spec}(k)$. In this case the schemes $X$ and
$X'$ are isomorphic as abstract schemes (but not as $k$-schemes).
By slightly abusing the notation, we 
will skip the subscript at the absolute Frobenius morphism and denote
it simply ${\sf F}$, as the relative Frobenius morphism. More
generally, for any $m\geq 1$ one defines $m$-th Frobenius twists $X^{(m)}$
and there is a morphism ${\sf F}^m\colon X\rightarrow X^{(m)}$, where ${\sf
  F}^m = {\sf F}\circ \dots \circ {\sf F}$ ($m$ times).
Let $\omega _X$ be the canonical invertible sheaf on $X$. 
Recall that the duality theory for finite flat
morphisms \cite{Har} yields that a right adjoint
functor ${{\sf F}^m}^{!}$ to ${\sf F}^m_{\ast}$ is isomorphic to
\begin{equation}\label{eq:rightadjointequat}
{{\sf F}^m}^{!}(?) = {{\sf F}^{m}}^{\ast}(?)\otimes \omega _{X/X^{(m)}} =
{{\sf F}^m}^{\ast}(?)\otimes \omega _{X}^{1-p^{m}}.
\end{equation}

\subsection{Koszul resolutions}\label{subsec:Kosres}

We recall here some basics of linear algebra. Let $\sf V$ be a finite dimensional vector
space over $k$ with a basis $\{e_1,\dots ,e_n\}$. Recall that the {\it $r$-th exterior power} $\wedge
^{r} \sf V$ of $\sf V$ is defined to be the $r$-th tensor power ${\sf V}^{\otimes
  r}$ of $\sf V$ divided by the vector subspace spanned by the elements:
$$ u_1\otimes \dots \otimes u_r - (-1)^{{\rm sgn} \sigma}u_{\sigma
  _1}\otimes \dots \otimes u_{\sigma (r)}$$

for all the permutations $\sigma \in \Sigma _r$ and
  $u_1,\dots ,u_r \in \sf V$. Similarly, the {\it $r$-th symmetric power}
  ${\sf S}^{r}\sf \V$ of $\sf V$ is defined to be the $r$-th tensor power
  $\V^{\otimes r}$ of $\sf V$ divided by the vector subspace spanned by the
  elements 
$$ u_1\otimes \dots \otimes u_r - u_{\sigma
  _1}\otimes \dots \otimes u_{\sigma (r)}$$

for all the permutations $\sigma \in \Sigma _r$ and
  $u_1,\dots ,u_r \in \sf V$. Finally, the {\it $r$-th divided power} ${\sf D}^r\sf V$ of $\sf V$ is
defined to be the dual of the symmetric power:
$$
{\sf D}^r\sf V = ({\sf S}^r{\sf V}^{\ast})^{\ast}.
$$

Let 
\begin{equation}\label{eq:seqvecsp}
0\rightarrow {\sf V}'\rightarrow \V \rightarrow {\sf V}''\rightarrow 0
\end{equation}

be a short exact sequence of vector spaces. For any $n > 0$
there is a functorial exact sequence (the Koszul resolution,
(\cite{Jan}, II.12.12))

\begin{equation}\label{eq:Kosqresol}
\dots \rightarrow {\sf S}^{n-i}{\V}\otimes \wedge ^i{\V}'\rightarrow \dots
\rightarrow {\sf S}^{n-1}{\V}\otimes {\V}'\rightarrow {\sf S}^n{\V}\rightarrow
{\sf S}^n{\V}''\rightarrow 0.
\end{equation}

Another fact about symmetric and exterior powers is the following (\cite{Har}, Exercise
5.16). For a short exact sequence as (\ref{eq:seqvecsp}) one has for each $n$ the
filtrations

\begin{equation}
{\sf S}^{n}{\V} = {\rm F}_n\supset {\rm F}_{n-1}\supset \dots \quad \mbox{and} \quad \bigwedge
^{n}{\V} ={\rm F}'_n\supset {\rm F}'_{n-1}\supset \dots
\end{equation}

such that 

\begin{equation}
{\rm F}_i/{\rm F}_{i-1}\simeq {\sf S}^{n-i}{\V}'\otimes {\sf S}^{i}{\V}''
\end{equation}

and 
\begin{equation}
{\rm F}'_i/{\rm F}'_{i-1}\simeq \bigwedge ^{n-i}{\sf V}'\otimes
\bigwedge ^{i}{\sf V}''
\end{equation}

When either $\V'$ or $\V''$ is a one-dimensional vector space, these
filtrations on exterior powers of $V$ degenerate into short exact sequences. 
If $\V''$ is one-dimensional, then one obtains:

\begin{equation}
0 \rightarrow \wedge ^{r}{\V}'\rightarrow \wedge ^{r}{\V} \rightarrow \wedge
^{r-1}{\V}'\otimes {\V}''\rightarrow 0.
\end{equation}

Similarly, if $\V'$ is one-dimensional, the filtration above
degenerates to give a short exact sequence:

\begin{equation}
0 \rightarrow \wedge ^{r-1}{\V}''\otimes {\V}'\rightarrow \wedge ^{r}{\V} \rightarrow \wedge
^{r}{\V}''\rightarrow 0.
\end{equation}

\subsection{Differential operators (\cite{Gro}, \cite{Haas})}\label{subsec:diffoper}
\subsubsection{True differential operators}
Let $X$ be a smooth scheme over $k$. Consider the product $X\times X$ and the diagonal $\Delta \subset
X\times X$. Let ${\mathcal J}_{\Delta}$ be the sheaf of ideals of
$\Delta$. 
\begin{definition}
An element $\phi \in {\mathcal End}_{k}(\Oo _X)$ is called a
differential operator if there exists some integer $n\geq 0$ such that
\begin{equation}
{\mathcal J}_{\Delta}^{n}\cdot \phi = 0.
\end{equation}
\end{definition}

One obtains a sheaf $\D _X$, the sheaf of differential
operators on $X$. Denote ${\mathcal
  J}_{\Delta}^{(n)}$ the sheaf of ideals generated by elements 
$a^n$, where $a\in {\mathcal J}_{\Delta}$. There is a filtration
on the sheaf $\D _X$ given by
\begin{equation}
\D _X^{(n)} = \{\phi \in {\mathcal End}_{k}(\Oo _X)\colon {\mathcal
  J}_{\Delta}^{(n)}\cdot \phi = 0\}.
\end{equation}

Since $k$ has characteristic $p$, one checks:
\begin{equation}\label{eq:p-filtration}
\D _X^{(p^{n})} = {\mathcal End}_{\Oo _X^{p^n}}(\Oo _X).
\end{equation}

Indeed, the sheaf ${\mathcal J}_{\Delta}$ is
generated by elements $a\otimes 1 - 1\otimes a$, where $a\in \Oo _X$, hence the sheaf ${\mathcal J}_{\Delta}^{(p^{n})}$ is
  generated by elements $a^{p^{n}}\otimes 1 - 1\otimes
  a^{p^{n}}$. This implies (\ref{eq:p-filtration}). One also checks that this
  filtration exhausts the whole $\D _X$, so one has (Theorem 1.2.4, \cite{Haas}):
\begin{equation}\label{eq:p-filtration}
\D _X = \bigcup _{n\geq 1}  {\mathcal End}_{\Oo _X^{p^n}}(\Oo _X).
\end{equation}

The filtration from (\ref{eq:p-filtration}) was called the
$p$--filtration in {\it loc.cit}. By definition of the Frobenius morphism
one has ${\rm H}^{i}(X,{\mathcal End}_{\Oo _X^{p^n}}(\Oo _X)) = {\rm
  H}^{i}(X^{(n)},{\mathcal End}({\sf F}^n_{\ast}\Oo _X))$. The sheaf $\D _X$
contains divided powers of vector fields (hence the name ``true'') as
opposed to the sheaf of PD-differential operators $\Dd _X$ that is discussed
in the next section.

\subsubsection{Crystalline differential operators}\label{subsec:Berthdiffoper}
The material of this subsection is taken from \cite{BMR}.
We recall, following {\it loc.cit.} the basic properties of crystalline differential operators (differential
operators without divided powers, or PD-differential
operators in the terminology of Berthelot and Ogus).\\

Let $X$ be a smooth variety, ${\T}^{\ast}_X$ the cotangent bundle, and 
${\rm T}^{\ast}(X)$ the total space of ${\T}^{\ast}_{X}$.

\begin{definition}
The sheaf $\Dd _X$ of {\it crystalline differential operators}
on $X$ is defined as the enveloping algebra of the
tangent Lie algebroid, i.e.,  for an affine open $U\subset X$ the
algebra $\Dd (U)$ contains the subalgebra $\Oo $ of functions, has an
$\Oo$-submodule identified with the Lie algebra of vector fields
$Vect(U)$ on $U$, and these subspaces generate $\Dd (U)$ subject to
relations $\xi_1\xi_2-\xi_2\xi_1=[\xi_1,\xi_2]\in Vect(U)$ for
$\xi_1, \xi_2\in Vect(U)$, and  $\xi \cdot f -f\cdot \xi =\xi(f)$
for $\xi\in Vect(U)$ and $f\in \Oo (U)$.
\end{definition}

Let us list the basic properties of the sheaf $\Dd _X$ \cite{BMR}:

\begin{itemize}

\item The sheaf of non-commutative algebras ${\sf F}_{\ast}\Dd _X$ has
  a center, which is isomorphic to $\Oo _{{\rm T}^{\ast}(X')}$, the
  sheaf of functions on the cotangent bundle to the Frobenius twist of
  $X$. The sheaf ${\sf F}_{\ast}\Dd _X$ is finite over its center.

\item This makes ${\sf F}_{\ast}\Dd _X$ a coherent sheaf on ${\rm
    T}^{\ast}(X')$. Thus, there exists a sheaf of algebras $\DD
 _X$ on ${\rm T}^{\ast}(X')$ such that $\pi _{\ast}\DD _X = {\sf
 F}_{\ast}{\rm D} _X$  (by abuse of notation we denote the projection ${\rm T}^{\ast}(X')\rightarrow X'$
 by the same letter $\pi$). The sheaf $\DD _X$ is an Azumaya algebra
    over ${\rm T}^{\ast}(X')$ of rank  $p^{2{\rm dim}(X)}$. 

\item There is a filtration on the sheaf ${\sf F}_{\ast}{\rm D}_X$
  such that the associated graded ring ${\rm gr}({\sf F}_{\ast}{\rm
  D}_X)$ is isomorphic to ${\sf F}_{\ast}\pi _{\ast}\Oo
 _{{\rm T}^{\ast}(X)} = {\sf F}_{\ast}{\sf S}^{\bullet}\T _X$.

\item Let $i\colon X'\hookrightarrow {\rm T}^{\ast}(X')$ be the zero section embedding.
Then $i^{\ast}\DD _X$ splits
 as an Azumaya algebra, the splitting bundle being ${\sf F}_{\ast}\Oo
 _X$. In other words, $i^{\ast}\DD _X = {\mathcal End}({\sf
 F}_{\ast}\Oo _X)$.        

\end{itemize}

Finally, recall that the sheaf $\Dd _X$ acts on $\Oo _X$ and that this action is not
faithful. It gives rise to a map $\Dd _X\rightarrow \D _X$;
its image is the sheaf of ``small differential operators'' ${\mathcal
  End}_{\Oo _X^p}(\Oo _X)$.

\subsection{Vanishing theorems for line bundles}\label{subsec:vanishingtheorems}
Let $\bf G$ be a connected, simply connected, semisimple algebraic
group over $k$, $\bf B$ a Borel subgroup of $\bf G$, and $\bf T$ a maximal
torus. Let $R({\bf T},{\bf G})$ be the root system of $\bf G$ with
respect to $\bf T$, $R^{+}$ the subset of positive roots, $S\subset
R^{+}$ the simple roots, and $h$ the Coxeter number of $\bf G$ that is
equal to $\sum m_i$, where $m_i$ are the coefficients of the highest
root of $\bf G$ written in terms of the simple roots $\alpha _i$.
 By $\langle \cdot, \cdot \rangle$ we denote the natural pairing $X({\bf
  T})\times Y({\bf T})\rightarrow {\mathbb Z}$, where $X({\bf T})$ is
the group of characters (also identified with the weight lattice) and
$Y({\bf T})$ the group of one parameter subgroups of $\bf T$ (also
identified with the coroot lattice). For a
subset $I\subset S$ let ${\bf P} = {\bf P}_{I}$ denote the associated
parabolic subgroup. Recall that the group of characters $X({\bf P})$
of $\bf P$ can be identified with $\{\lambda \in X({\bf T})|\langle
\lambda ,\alpha ^{\vee}\rangle = 0,$ for all $\alpha \in I\}$. In
particular, $X({\bf B}) = X({\bf T})$. A weight $\lambda \in X({\bf
  B})$ is called {\it dominant} if $\langle \lambda ,\alpha
^{\vee}\rangle\geq 0$ for all $\alpha \in S$. A dominant weight $\lambda \in
X({\bf P})$ is called $\bf P$-regular if $\langle \lambda ,\alpha
^{\vee}\rangle> 0$ for all $\alpha \notin I$, where ${\bf P} = {\bf P}_I$ is
a parabolic subgroup. A weight $\lambda$ defines a line bundle $\Ll
_{\lambda}$ on ${\bf G}/{\bf B}$. Line bundles on ${\bf G}/{\bf B}$
that correspond to dominant weights are ample. If a weight $\lambda$ is ${\bf
P}$-regular, then the corresponding line bundle is ample on ${\bf G}/{\bf P}$.
Finally, denote $\mathcal W$ the Weyl group of the root system $R({\bf T},{\bf G})$.

\subsubsection{Cohomology of line bundles on flag varieties}

Recall first the Borel--Weil--Bott theorem in the characteristic zero case. For $w\in
{\mathcal W}$ denote $l(w)$ the length of $w$. Let $\chi$ be a
weight. For an element $w\in {\mathcal W}$ define the dot
action of $w$ on a weight $\chi \in X({\bf T})$ by the rule $w\cdot
\chi =w(\chi + \rho) - \rho$. If $\chi + \rho \in X_{+}({\bf T})$ and the ground field $k$
has characteristic zero, then the Borel--Weil--Bott theorem asserts that ${\rm H}^i({\bf G}/{\bf
  B},\Ll _{w\cdot \chi})=0$ for $i\neq l(w)$ and ${\rm H}^{l(w)}({\bf G}/{\bf
  B},\Ll _{w\cdot \chi})={\rm H}^0({\bf G}/{\bf B},\Ll _{\chi})$. Let
now the field $k$ is of characteristic $p>0$. Recall that the dominant
bottom $p$--alcove is defined to be
\begin{equation}
{\rm C} = \{\lambda \in X({\bf T}) \ | \ 0< \langle \lambda + \rho,\alpha
^{\vee}\rangle < p, \forall \ \alpha \in {\rm R}_{+}\},
\end{equation}

and its closure 
\begin{equation}
{\rm \bar C} = \{\lambda \in X({\bf T}) \ | \ 0\leq  \langle \lambda + \rho,\alpha
^{\vee}\rangle \leq p, \forall \ \alpha \in {\rm R}_{+}\}.
\end{equation}

First, one has (\cite{And}, Theorem 2.3):

\begin{theorem}\label{th:AndthII}
  If $\chi \in {\rm \bar C}$, then 
\begin{equation}
{\rm H}^i({\bf G}/{\bf B},\Ll _{\chi}) = {\rm H}^{i+1}({\bf G}/{\bf
  B},\Ll _{s_{\alpha}\cdot \chi}).
\end{equation}
\end{theorem}

More generally, 

\begin{theorem} {\rm(}\cite{And}{\rm )}
If $\chi \in {\rm \bar C}$ and $w\in \mathcal W$, then ${\rm H}^{l(w)}({\bf G}/{\bf
  B},\Ll _{w\cdot \chi})={\rm L}_{\lambda}$ and zero otherwise. Here
  ${\rm L}_{\lambda}$ is the simple $\bf G$-module with highest weight
  $\chi$.
\end{theorem}

\begin{definition}
A weight $\chi \in X({\bf T})$ is called generic if the Borel--Weil--Bott theorem holds
for $\Ll _{\chi}$, that is if $j=l(w)$ is the sole degree of
cohomology for which ${\rm H}^j({\bf G}/{\bf B},\Ll _{\chi})$ does not vanish.
\end{definition}

\subsubsection{First cohomology group of a line bundle}

Consider a simple root $\alpha \in {\rm S}$. It defines a reflection $s_{\alpha}$ in $\mathcal W$. One has $s_{\alpha}\cdot
\chi = s_{\alpha}(\chi)-\alpha$. Recall Andersen's theorem \cite{An}
on the first cohomology group of a line bundle on ${\bf G}/{\bf B}$. Let
$\chi \in X({\bf T})$ be a weight and $\Ll _{\chi}$ the corresponding
line bundle on ${\bf G}/{\bf B}$.

\begin{theorem}\label{th:Andth}
The group ${\rm H}^1({\bf G}/{\bf B},\Ll _{\chi})$ is non--vanishing if and only if
there exist a simple root $\alpha$ such that one of the following conditions is satisfied:

\begin{itemize}

\item $-p\leq \langle \chi ,\alpha ^{\vee}\rangle \leq -2$ and
  $s_{\alpha}\cdot \chi = s_{\alpha}(\chi)-\alpha$ is dominant.

\item $\langle \chi ,\alpha ^{\vee}\rangle = -ap^n-1$ for some $a,n\in
  {\bf N}$ with $a<p$ and $s_{\alpha}(\chi)-\alpha$ is dominant.

\item $-(a+1)p^n\leq  \langle \chi ,\alpha ^{\vee}\rangle \leq
  -ap^n-2$ for some $a,n\in {\bf N}$ with $a<p$ and $\chi +
  ap^n\alpha$ is dominant.
\end{itemize}
\end{theorem}

Further, Corollary 3.2 of \cite{An} states:
\begin{theorem}\label{th:Andcor}
Let $\chi$ be a weight. If either $\langle \chi ,\alpha ^{\vee}\rangle
\geq -p$ or  $\langle \chi ,\alpha ^{\vee}\rangle = -ap^n-1$ for some
$a, n\in {\bf N}$ and $a<p$, then 
\begin{equation}
{\rm H}^i({\bf G}/{\bf B},\Ll _{\chi}) = {\rm H}^{i-1}({\bf G}/{\bf
  B},\Ll _{s_{\alpha}\cdot \chi}).
\end{equation}
\end{theorem}

\subsubsection{Line bundles on the cotangent bundles of flag varieties}

In this section we follow \cite{KLT}. Recall that the prime $p$ is {\it a good prime} for $\bf G$ if $p$ is
coprime to all the coefficients of the highest root of $\bf G$ written
in terms of the simple roots. In particular, if $\bf G$ is a simple
group of type $\bf A$, then all primes are good for $\bf G$; if
$\bf G$ is either of the type $\bf B$ or $\bf D$, then good primes are
$p\geq 3$; if $\bf G$ is of type ${\bf G}_2$, then $p$ is a good prime
for $\bf G$ if $p\geq 6$.\\

Recall first the Kempf vanishing theorem \cite{Ke}:
\begin{theorem}
Let $\lambda \in X({\bf T})$ be a dominant weight, and $\Ll
_{\lambda}$ the associated line bundle on ${\bf G}/{\bf B}$. Then
\begin{equation}
{\rm H}^i({\bf G}/{\bf B},\Ll _{\lambda}) = 0
\end{equation}

for $i>0$.

\end{theorem}

The next theorem concerns with the vanishing of cohomology groups of line
bundles on the total spaces of cotangent bundles (\cite{KLT}, Theorem 5):

\begin{theorem}\label{th:KLTvantheorem}
Let ${\rm T}^{\ast}({\bf G}/{\bf B})$ be the total space of the
  cotangent bundle of the flag variety ${\bf G}/{\bf B}$, and
  $\pi \colon {\rm T}^{\ast}({\bf G}/{\bf B})\rightarrow {\bf G}/{\bf B}$ the projection.
Assume that $\mbox{char} \ k$ is a good prime for $\bf G$. Let $\lambda\in
X({\bf B})$ be a weight such that $\langle \lambda ,\alpha
  ^{\vee}\rangle \geq -1$ for $\forall \alpha \in {\rm R}^{+}$. Then 
\begin{equation}\label{eq:GmodPvanishing}
{\rm H}^i({\rm T}^*({\bf G}/{\bf B}),\Ll (\lambda)) = 
{\rm H}^{i}({\bf G}/{\bf B},\Ll _{\lambda}\otimes {\sf
  S}^{\bullet}{\mathcal T}_{{\bf G}/{\bf B}}) = 0
\end{equation}

\noindent for $i>0$.
\end{theorem}

\noindent In particular, one has:
\begin{equation}\label{eq:OacyclicGmodB}
{\rm H}^i({\rm T}^*({\bf G}/{\bf B}),\Oo _{{\rm T}^*({\bf G}/{\bf B})}) = {\rm H}^{i}({\bf G}/{\bf B},{\sf
  S}^{\bullet}{\mathcal T}_{{\bf G}/{\bf B}}) = 0
\end{equation}

\noindent for $i>0$.\\

The tangent bundle ${\mathcal T}_{{\bf G}/{\bf B}}$ is a
  homogeneous bundle associated with the $\bf B$-module ${\mathfrak
  g}/{\mathfrak b}$ under the adjoint action. There is an isomorphism $({\mathfrak g}/{\mathfrak
  b})^{\ast}\simeq {\mathfrak u}$ of $\bf B$-modules in good
  characteristic. The following lemma is a straightforward modification
  of an argument used in \cite{KLT} in the proof of Theorem \ref
  {th:KLTvantheorem}); for convenience of the reader we give the
  proof. To simplify the notation, we
  will drop the superscript that denotes the Frobenius twist of a
  scheme, and will use the absolute Frobenius morphism (it does not
  affect the cohomology groups).

\begin{lemma}\label{lem:KLTlemma}
Let $\lambda$ be a weight and $\Ll _{\lambda}$ is the line bundle
associated to $\lambda$. Then for any fixed $i>0$ one has:
\begin{equation}
{\rm H}^{i+j}({\bf G}/{\bf B},{\sf F}_{\ast}\Omega _{{\bf G}/{\bf B}}^j\otimes \Ll
_{\lambda}) = 0 \ \mbox{\rm for all} \ j\geq 0 \ \Rightarrow {\rm
  H}^{i}({\bf G}/{\bf B},{\sf F}_{\ast}{\sf
  S}^{\bullet}{\mathcal T}_{{\bf G}/{\bf B}}\otimes \Ll _{\lambda}) = 0
\end{equation}
\end{lemma}

\begin{proof}
Consider the short exact sequence of $\bf B$-modules:
\begin{equation}
0\rightarrow ({\mathfrak b}/{\mathfrak u})^{\ast}\rightarrow
{\mathfrak b}^{\ast}\rightarrow {\mathfrak u}^{\ast}\rightarrow 0.
\end{equation}

The $\bf B$-module $({\mathfrak b}/{\mathfrak u})^{\ast}$ is
trivial, hence induces a trivial bundle on ${\bf G}/{\bf B}$.  The Koszul resolution associated to this sequence reads:
\begin{equation}
\dots \rightarrow \wedge ^{j}({\mathfrak b}/{\mathfrak
  u})^{\ast}\otimes {\sf S}^{k-j}{\mathfrak b}^{\ast}\rightarrow \dots
  \rightarrow ({\mathfrak b}/{\mathfrak u})^{\ast}\otimes {\sf
  S}^{k-1}{\mathfrak b}^{\ast}\rightarrow {\sf S}^{k}{\mathfrak
  b}^{\ast}\rightarrow {\sf S}^{k}{\mathfrak u}^{\ast}\rightarrow 0.
\end{equation}

Apply ${\sf F}_{\ast}$ to this resolution and tensor the obtained
complex with $\Ll _{\lambda}$. For a fixed $l\geq 0$ the vanishing of
the groups ${\rm H}^{i}({\bf G}/{\bf B},{\sf F}_{\ast}{\sf
  S}^{k}{\mathfrak b}^{\ast}\otimes \Ll _{\lambda})$ for all $k\geq 0$
  and $i>l$ implies the vanishing of the group ${\rm H}^{i}({\bf
    G}/{\bf B},{\sf F}_{\ast}{\sf
  S}^{k}{\mathfrak u}^{\ast}\otimes \Ll _{\lambda})$ for all $k\geq 0$
  and $i>l$. Consider now the short exact sequence
\begin{equation}
0\rightarrow ({\mathfrak g}/{\mathfrak b})^{\ast}\rightarrow
{\mathfrak g}^{\ast}\rightarrow {\mathfrak b}^{\ast}\rightarrow 0,
\end{equation}

and associated Koszul complex:
\begin{equation}
\dots \rightarrow \wedge ^{j}({\mathfrak g}/{\mathfrak
  b})^{\ast}\otimes {\sf S}^{k-j}{\mathfrak g}^{\ast}\rightarrow \dots \rightarrow ({\mathfrak g}/{\mathfrak
  b})^{\ast}\otimes {\sf S}^{k-1}{\mathfrak g}^{\ast}\rightarrow {\sf S}^k{\mathfrak g}^{\ast}\rightarrow {\sf
  S}^{k}{\mathfrak b}^{\ast}\rightarrow 0.
\end{equation}

Again apply ${\sf F}_{\ast}$ to this resolution and tensor the obtained
complex with $\Ll _{\lambda}$.
The $\bf B$-module $\wedge ^j({\mathfrak g}/{\mathfrak b})^{\ast}$
induces the bundle of $j$-forms $\Omega _{{\bf G}/{\bf B}}^j$ on ${\bf
  G}/{\bf B}$. Therefore, for a fixed $i$ we have an implication:
\begin{equation}
{\rm H}^{i+j}({\bf G}/{\bf B},{\sf F}_{\ast}\Omega _{{\bf G}/{\bf B}}^j\otimes \Ll
_{\lambda}) = 0 \ \mbox{\rm for all} \ j\geq 0 \ \Rightarrow {\rm
  H}^{i}({\bf G}/{\bf B},{\sf F}_{\ast}{\sf
  S}^{\bullet}{\mathfrak b}^{\ast}\otimes \Ll _{\lambda}) = 0,
\end{equation}

which in turn implies ${\rm H}^{i}({\bf G}/{\bf B},{\sf F}_{\ast}{\sf
  S}^{\bullet}{\mathfrak u}^{\ast}\otimes \Ll _{\lambda}) = 0$.
\end{proof}

\subsection{Frobenius splitting and {\rm F}-amplitude (\cite{Ar}, \cite{Smith})}\label{subsec:Frobsplit}

We first recall several well-known definitions concerning Frobenius
splitting and $\rm F$-regularity \cite{Smith}.

\begin{definition}
A variety $X$ is called Frobenius split if the natural map $\Oo
_X\rightarrow {\sf F}_{\ast}\Oo _X$ is split (as a map of $\Oo _X$-modules).
\end{definition}

We will need a refinement of the previous definition, which is the
Frobenius splitting along a divisor. Let $D$ be an effective
Cartier divisor, and $s$ be a section of $\Oo _X(D)$ vanishing
precisely along $D$. Consider the map $\Oo _X\rightarrow \Oo _X(D)$
sending $1\rightarrow s$. This induces a map of $\Oo _X$-modules
\begin{equation}
\Oo _X\rightarrow {\sf F}_{\ast}\Oo _X\rightarrow {\sf F}_{\ast}\Oo _X(D),
\end{equation}

where the first arrow is the Frobenius map and the second arrow is (the pushforward
of ) the map $1\rightarrow s$. Denote by $^1s$ the element $s$
considered as an element of ${\sf F}_{\ast}\Oo _X(D)$. The variety $X$ is said to be Frobenius $D$-split if this composition map splits, that is,
if there is a map ${\sf F}_{\ast}\Oo _X(D)\rightarrow \Oo _X$ sending
$^1s\rightarrow 1$.

There is a stronger notion of $\rm F$-regular varieties. One of the
equivalent formulations is that $X$ is $\rm F$-regular if it is stably
Frobenius split along every effective Cartier divisor (\cite{Smith},
Theorem 3.10). In {\it loc.cit.} it is also proved that smooth Fano
varieties are $\rm F$-regular (for $p\gg 0$).

\begin{definition}\label{def:f-amplitude}
The $\rm F$-amplitude $\phi (\Ee)$ of a coherent
  sheaf $\Ee$ is the smallest integer $l$ such that, for any locally
  free sheaf $\Ff$, there exists an integer $N$ such that 
  ${\rm H}^{i}(X,{\sf F}_{m}^{\ast}(\Ee)\otimes \Ff) = 0$ for $i>l$
  and $m > N$.
\end{definition}

A sheaf $\Ee$ is called {\it $\rm F$-ample} if
its $\rm F$-amplitude $\phi (\Ee)$ is equal to zero. 

\begin{proposition} (\cite{Ar}, Proposition 8.1)\label{lem:strongervanishing}
Let $X$ be a Frobenius split variety, and let $\mathcal G$ be a
  coherent sheaf on $X$. Then

\begin{itemize}

\item [(i)] ${\rm H}^{i}(X,{\mathcal G}) = 0$ for $i>\phi ({\mathcal G})$.

\item [(ii)] If ${\mathcal G}$ is locally free, then ${\rm H}^{i}(X,{\mathcal G}\otimes \omega _X) = 0$ for $i>\phi ({\mathcal G})$.

\end{itemize}
\end{proposition}

For convenience of the reader let us give the proof.

\begin{proof}

(i): For any $m\geq 0$ one has embeddings of cohomology groups:
\begin{equation}\label{eq:frobembed2}
{\rm H}^{k}(X,{\mathcal G})\hookrightarrow {\rm H}^{k}(X,{\sf F}_{m}^{\ast}{\mathcal G})
\end{equation}

\noindent By definition of the $\rm F$-amplitude the right-hand
side group is equal to zero for $k > \phi ({\mathcal G})$ and $m$ large enough, hence the statement
of (i).\\

\noindent (ii):  If $\mathcal G$ is locally free, then by Serre duality one has:
\begin{equation}\label{eq:ArSerdual}
{\rm H}^{l}(X,{\mathcal G}\otimes \omega _X) = {\rm H}^{n-l}(X,{\mathcal G}^{\ast})^{\ast}
\end{equation}

\noindent Taking the dual of (\ref{eq:frobembed2}) for ${\mathcal
  G}^{\ast}$ and $k = n-l$, for any $m$ we get a surjection
\begin{equation}\label{eq:Arsurjection}
{\rm H}^{n-l}(X,{\sf F}_{m}^{\ast}{\mathcal G}^{\ast})^{\ast} \twoheadrightarrow
{\rm H}^{n-l}(X,{\mathcal G}^{\ast})^{\ast}
\end{equation}

\noindent Applying again Serre duality to the right-hand side group in
(\ref{eq:Arsurjection}) and taking into account (\ref{eq:ArSerdual}) we get a
surjection
\begin{equation}
{\rm H}^{l}(X,{\sf F}_{m}^{\ast}{\mathcal G}\otimes \omega
_X)\twoheadrightarrow {\rm H}^{l}(X,{\mathcal G}\otimes \omega _X)
\end{equation}

\noindent By definition of the $\rm F$-amplitude the left-hand
side group is equal to zero for $l > \phi ({\mathcal G})$ and $m$ large enough, hence the statement
of (ii).

\end{proof}

Assume $X$ is Frobenius split.

\begin{proposition}\label{lem:vanishviaomega}
Let $l = \phi ({{\sf F}_m}_{\ast}{\omega _{X}^{-p^m}})$. Then $\Ext
  ^{k}({{\sf F}_m}_{\ast}{\Oo _{X}},{{\sf F}_m}_{\ast}{\Oo _{X}}) = 0$
  for $k>l$.
\end{proposition}

\begin{proof}
Recall that according to (\ref{eq:ext-cohom}) one has 
\begin{equation}\label{eq:exttocohom1}
\Ext ^{k}_{X}({{\sf F}_m}_{\ast}\Oo _{X},{{\sf F}_m}_{\ast}\Oo _{X}) = 
{\rm H}^{k}(X,{\sf F}_m^{\ast}{{\sf F}_m}_{\ast}{\Oo _{X}}\otimes \omega _{X}^{1-p^m})
\end{equation}

The isomorphism of sheaves ${\sf F}_m^{\ast}{{\sf F}_m}_{\ast}{\Oo _{X}}\otimes
\omega _{X}^{1-p^m} = {\sf F}_m^{\ast}({{\sf F}_m}_{\ast}\Oo _{X}\otimes \omega
_{X}^{-1})\otimes \omega _{X} = {\sf F}_m^{\ast}{{\sf F}_m}_{\ast}\omega
_{X}^{-p^m}\otimes \omega _{X}$ implies
\begin{equation}\label{eq:exttocohom2}
\Ext ^{k}_{X}({{\sf F}_m}_{\ast}\Oo _{X},{{\sf F}_m}_{\ast}\Oo _{X}) =
{\rm H}^{k}(X,{\sf F}_m^{\ast}{{\sf F}_m}_{\ast}\omega
_{X}^{-p^m}\otimes \omega _{X})
\end{equation}

Note that for a coherent sheaf $\Ee$ on $X$ and $n\geq 1$ one has $\phi ({\sf
  F}_n^{\ast}\Ee) = \phi (\Ee)$. Indeed, denote $\phi _n = \phi ({\sf
  F}_n^{\ast}\Ee)$. If $k>\phi _n$ and $\Ff$ is a coherent sheaf,
  then, by definition, for sufficiently large $l$ one has ${\rm H}^{k}(X,{\sf
  F}_l^{\ast}{\sf F}_n^{\ast}\Ee\otimes \Ff)=  {\rm H}^{k}(X,{\sf
  F}_{l+n}^{\ast}\Ee\otimes \Ff)=0$. Thus, $\phi _n\leq
  \phi$. The inverse inequality is similar. Now by Lemma \ref{lem:strongervanishing} for $k > \phi ({\sf
  F}_m^{\ast}{{\sf F}_m}_{\ast}{\omega _{X}^{-p^m}}) = \phi
  ({{\sf F}_m}_{\ast}{\omega _{X}^{-p^m}})$ the right-hand side in (\ref{eq:exttocohom2}) is equal to zero,
  hence the statement.

\end{proof}

In particular, if $\phi ({{\sf F}_m}_{\ast}{\omega
  _{X}^{-p^m}})=\phi ({{\sf F}_m}_{\ast}\Oo _X\otimes \omega _X^{-1}) =
  0$, then $\Ext ^{k}_{X}({{\sf F}_m}_{\ast}\Oo
  _{X},{{\sf F}_m}_{\ast}\Oo _{X}) = 0$ for $k>0$.
In general, the ${\rm F}$--ampleness is too strong a property for the
  bundle ${{\sf F}_m}_{\ast}\Oo _X\otimes \omega _X^{-1}$ 
(for example, the $\rm F$-amplitude of ${{\sf F}_m}_{\ast}\Oo _X\otimes
\omega _X^{-1}$ is zero in the simplest example of
$\Pp ^n$ and greater than zero already in the case of 3-dimensional
quadrics ${\sf Q}_3$, while higher $\Ext$-groups from
  (\ref{eq:exttocohom2}) vanish in this case \cite{Sam1},
  \cite{SamD-affflag}). However, if the variety has
additional properties, apart from being Frobenius split (e.g., being a
toric variety), then one can
derive some vanishing statements from Proposition \ref{lem:vanishviaomega}
(see Section \ref{sec:toricFanovar}).

\subsection{Derived categories of coherent sheaves (\cite{Or})}\label{subsec:dercatofcohsheaves}

In this section we recall some facts about semiorthogonal
decompositions in derived categories of coherent sheaves and tilting
equivalences. We refer the reader to {\it loc.cit.} for the definition of
semiorthogonal decompositions in derived categories.

\subsubsection{Semiorthogonal decompositions}\label{subsubsec:semdec}

Let $S$ be a smooth scheme, and $\Ee$ a vector bundle of rank $n$ over $S$. Denote $X = \Pp _{S}(\Ee)$ the projectivization of the
bundle $\Ee$. Let $\pi \colon X\rightarrow S$ be the projection, and
$\Oo _{\pi}(-1)$ the relative invertible sheaf on $X$.
\begin{theorem}\label{th:Orlovth1}
The category $\Dd ^{b}(X)$ admits a semiorthogonal decomposition:
\begin{equation}
\Dd ^{b}(X) = \langle \pi ^{\ast}\Dd ^{b}(S)\otimes \Oo
_{\pi}(-n+1),\pi ^{\ast}\Dd ^{b}(S)\otimes \Oo _{\pi}(-n+2),\dots,\pi
^{\ast}\Dd ^{b}(S)\rangle .
\end{equation}
\end{theorem}

Further, we need a particular case of another Orlov's
theorem ({\it loc.cit.}) Consider a smooth scheme $X$ and a closed smooth
subscheme $i\colon Y\subset X$ of codimension two. Let $\tilde X$ be the blowup of $X$
along $Y$. There is a cartesian square:

\begin{figure}[H]
$$
\xymatrix @C5pc @R4pc {
 {\tilde Y}\ar[r]^{j} \ar@<-0.1ex>[d]^{p} & {\tilde X} \ar@<-0.4ex>[d]^{\pi} \\
       Y   \ar[r]^{i} & X
 }
$$
\end{figure}

Here $\tilde Y$ is the exceptional divisor. If ${\mathcal
  N}_{Y/X}$ is normal bundle to $Y$ in $X$, then the projection
$p$ is the projectivization of the bundle ${\mathcal N}_{Y/X}$. Denote $\Oo
_{p}(-1)$ be the relative invertible sheaf with respect to $p$.

\begin{theorem}\label{th:Orlovth2}
The category $\Dd ^{b}(\tilde X)$ admits a semiorthogonal
  decomposition:
\begin{equation}
\Dd ^{b}(\tilde X) = \langle j_{\ast}(p^{\ast}\Dd ^{b}(X)\otimes \Oo
_{p}(-1)), \pi ^{\ast}\Dd ^{b}(X)\rangle .
\end{equation}
\end{theorem}

\subsubsection{Tilting equivalences}

\begin{definition}\label{def:tiltingdefinition}
A coherent sheaf $\Ee$ on a smooth variety $X$ over an algebraically
closed field $k$ is called a {\sf tilting generator} of
the bounded derived category $\D ^b(X)$ of coherent sheaves on
$X$ if the following holds:
\begin{enumerate}
\item The sheaf $\Ee$ is a tilting object in $\D^b(X)$, that
  is, for any $i \geq 1$ one has $\Ext^i(\Ee,\Ee)=0$.
\item The sheaf $\Ee$ generates the derived category $\D ^{-}(X)$
  of complexes bounded from above, that is, if for some object $\Ff
  \in \D^{-}(X)$ one has $\RHom ^{\bullet}(\Ff ,\Ee)=0$, then $\Ff=0$.
\end{enumerate}
\end{definition}

Tilting generators are a tool to construct derived equivalences. One has:

\begin{lemma}$\rm ($\cite{Kaledquant}, Lemma 1.2$\rm )$\label{lem:tiltinglemma}
Let $X$ be a smooth scheme, $\Ee$ a tilting generator of the derived category $\D^b(X)$, and
denote $R={\mathcal End}(\Ee)$. Then the
algebra $R$ is left-Noetherian, and the correspondence $\Ff \mapsto
\RHom^{\bullet}(\Ee,\Ff)$ extends to an equivalence
\begin{equation}\label{equi}
  \D^b(X) \rightarrow \D^b(R-\mbox{\rm mod}^{\rm fg})
\end{equation}
between the bounded derived category $\D^b(X)$ of coherent
sheaves on $X$ and the bounded derived category $\D^b(R-\mbox{\rm mod}^{\rm fg})$ of
finitely generated left $R$-modules.
\end{lemma}

One can drop the condition (2) in Definition \ref{def:tiltingdefinition} and only
require the vanishing of groups $\Ext^i(\Ee,\Ee)=0$ for $i>0$. Lemma
\ref{lem:tiltinglemma} claims in this case that there is a full
faithful embedding of the category  $\D^b(R\fmod)$ into $\D^b(X)$. 
If $\Ee$ is such a bundle, then we will call it {\it partial
  tilting}. The goal of this paper is to study when the bundle ${\sf
  F}_{\ast}\Oo _X$ gives a (partial) tilting bundle on $X$. 

Clearly, given a pair of smooth varieties $X$ and $Y$ such that both ${\sf
  F}_{\ast}\Oo _X$ and ${\sf F}_{\ast}\Oo _Y$ are (partial) tiltng
  bundles, the bundle ${\sf F}_{\ast}\Oo _{X\times Y}$ is also a (partial) tilting
  bundle on $X\times Y$. Indeed, $\Oo _{X\times Y}=\Oo _X\boxtimes
  _{k}\Oo _{Y}$. Therefore ${\sf F}_{\ast}\Oo _{X\times Y}={\sf
  F}_{\ast}\Oo _X\boxtimes _{k}{\sf F}_{\ast}\Oo _Y$. Hence,
\begin{equation}
\Ext ^i({\sf F}_{\ast}\Oo _{X\times Y},{\sf F}_{\ast}\Oo _{X\times
  Y})={\rm H}^i(X\times Y,(({\sf F}_{\ast}\Oo _X)^{\vee}\otimes {\sf
  F}_{\ast}\Oo _X)\boxtimes _{k}(({\sf F}_{\ast}\Oo _Y)^{\vee}\otimes {\sf
  F}_{\ast}\Oo _Y))=0 
\end{equation}

for $i>0$ by the K\"unneth formula. If both ${\sf F}_{\ast}\Oo _X$ and
${\sf F}_{\ast}\Oo _Y$ are generators in $\Dd ^b(X)$ and $\Dd ^b(Y)$,
respectively, then ${\sf F}_{\ast}\Oo _X\boxtimes _{k}{\sf
  F}_{\ast}\Oo _Y$ is a generator in $\Dd ^b(X\times Y)$. This is a
particular case of Lemma 3.4.1 in \cite{BonVdB}. 

\subsection{Derived localization theorem (\cite{BMR})}

We need to recall the derived Beilinson--Bernstein localization
theorem (\cite{BMR}). Let $\bf G$ be a semisimple algebraic group over $k$, ${\bf G}/{\bf
  B}$ the flag variety, and $\Uu ({\mathfrak g})$ the universal
  enveloping algebra of the corresponding Lie algebra. The center of
  $\Uu ({\mathfrak g})$ contains the ``Harish--Chandra part''
  ${\mathfrak Z}_{\rm HC}=\Uu ({\mathfrak g})^{\bf G}$. Denote $\Uu
  ({\mathfrak g})_0$ the central reduction $\Uu ({\mathfrak g})\otimes _{{\mathfrak
  Z}_{\rm HC}} {\bf k}$, where $\bf k$ is the trivial $\mathfrak g$-module.
Consider the category ${\rm D}_{{\bf G}/{\bf
B}}$-\mbox{\rm mod} of coherent ${\rm D}_{{\bf G}/{\bf B}}$-modules and the
category $\Uu ({\mathfrak g})_0$-\mbox{\rm mod} of finitely generated modules over $\Uu
  ({\mathfrak g})_0$. The
  derived localization theorem (Theorem 3.2, \cite{BMR}) states:

\begin{theorem}\label{th:derBB-equivalence}
Let char $k = p>h$, where $h$ is the Coxeter number of the group $\bf
G$. Then there is an equivalence of derived categories:
\begin{equation}\label{eq:derBB-equivalence}
\Dd ^{b}({\rm D}_{{\bf G}/{\bf B}}-\mbox{\rm mod}) \simeq \Dd ^{b}(\Uu ({\mathfrak g})_0-\mbox{\rm mod}),
\end{equation}

\end{theorem}

\subsection{Frobenius descent}
Let $\Ee$ be an $\Oo _X$-module equipped with an integrable connection
$\nabla$, and let $U\subset X$ be an open subset with local
coordinates $t_1,\dots ,t_d$. Let $\partial _1,\dots ,\partial _d$ be
the local basis of derivations that is dual to the basis $(dt_i)$ of
$\Omega ^1_X$. The connection $\nabla$ is said to have zero
$p$-curvature over $U$ if and only if for any local section $s$ of
$\Ee$ and any $i$ one has $\partial _i^ps=0$. For any $\Oo
_{X'}$-module $\Ee$ the Frobenius pullback ${\sf F}^{\ast}\Ee$ is
equipped with a canonical integrable connection with zero
$p$-curvature. The Cartier descent theorem (Theorem 5.1, \cite{Katz})
states that the functor ${\sf F}^{\ast}$ induces an equivalence
between the category of $\Oo _{X'}$-modules and that of $\Oo
_X$-modules equipped with an integrable connection with zero
$p$-curvature.

On the other hand, if an $\Oo _X$-module $\Ee$ is equipped with  an integrable connection with zero
$p$-curvature, then it has a structure of left $\Dd _X$-module. Given
the Cartier descent theorem, we see that the Frobenius pullback ${\sf
  F}^{\ast}\Ee$ of a coherent sheaf $\Ee$ on $X'$ is a left $\Dd
_X$-module.\\

There are also ``unbounded'' versions of Theorem
\ref{th:derBB-equivalence}, the bounded categories in
(\ref{eq:derBB-equivalence}) being replaced by categories unbounded from
above or below (see Remark 2 on p. 18 of \cite{BMR}). Thus extended, Theorem \ref{th:derBB-equivalence} implies:

\begin{lemma}\label{lem:Frobgeneratorforp>h}
Let $\bf G$ be a semisimple algebraic group over $k$, and $X={\bf
  G}/{\bf B}$ the flag variety. Let char $k = p >h$, where $h$ is the
  Coxeter number of the group ${\bf G}$. Then the bundle ${\sf F}_{\ast}\Oo _X$
satisfies the condition (2) of Definition \ref{def:tiltingdefinition}.
\begin{proof}
We \ need \ to \ show \ that \ for \ an \ object \ $\Ee
  \in \D^{-}(X')$ \ the \ equality \ $\RHom ^{\bullet}(\Ee ,{\sf F}_{\ast}\Oo
  _X)=0$ implies $\Ee=0$. By adjunction we get:
\begin{equation}
{\mathbb H}^{\bullet}(X,({\sf F}^{\ast}\Ee)^{\vee}) = 0.
\end{equation}

On the other hand,
\begin{eqnarray}
& ({\sf F}^{\ast}\Ee)^{\vee}={\mathcal R}{\mathcal
  Hom}_{\Oo _X}({\sf F}^{\ast}\Ee,\Oo _X) = {\mathcal R}{\mathcal
  Hom}_{\Oo _X}({\sf F}^{\ast}\Ee,{\sf F}^{\ast}\Oo _{X'})  = & \nonumber \\
& ={\sf F}^{\ast}{\mathcal R}{\mathcal Hom}_{\Oo _{X'}}(\Ee,\Oo _{X'})={\sf F}^{\ast}\Ee ^{\vee}.
\end{eqnarray}

By the Cartier descent, the object ${\sf F}^{\ast}\Ee ^{\vee}$ is an object of the category
$\Dd ^{-}(\Dd _X-\mbox{\rm mod})$. Now ${\sf F}^{\ast}\Ee ^{\vee}$ is annihilated by the functor
${\rm R}\Gamma$. Under our assumptions on $p$, this functor is an equivalence of categories by Theorem
\ref{th:derBB-equivalence}. Hence, ${\sf F}^{\ast}\Ee ^{\vee}$ is
quasiisomorphic to zero, and therefore so are $\Ee ^{\vee}$ and $\Ee$. 
\end{proof}
\end{lemma}

\begin{corollary}\label{cor:Frobgeneratorforp>h}
Let $\bf P$ be a parabolic subgroup of $\bf G$, and let $p>h$. Then
  the bundle ${\sf F}_{\ast}\Oo _{{\bf G}/{\bf P}}$ is a generator in ${\rm D}^{b}({\bf
  G}/{\bf P})$.
\end{corollary}

\begin{proof}
Denote $Y={\bf G}/{\bf P}$, and let $\pi : X = {\bf G}/{\bf B}\rightarrow Y$ be the projection. As before, one has to show
that for any object $\Ee \in \Dd ^{-}(Y')$ the equality $\RHom ^{\bullet}(\Ee ,{\sf F}_{\ast}\Oo
  _Y)=0$ implies $\Ee=0$. Notice that  ${\rm R}^{\bullet}\pi _{\ast}\Oo _{X} = \Oo _{Y}$. The condition $\RHom ^{\bullet}(\Ee ,{\sf F}_{\ast}\Oo
  _Y)=0$ then translates into:
\begin{equation}
\RHom ^{\bullet}(\Ee ,{\sf F}_{\ast}\Oo  _Y)= \RHom ^{\bullet}(\Ee,
{\sf F}_{\ast}{\rm R}^{\bullet}\pi _{\ast}\Oo _X) = \RHom ^{\bullet}(\Ee
,{\rm R}^{\bullet}\pi _{\ast}{\sf F}_{\ast}\Oo _X) =0.
\end{equation}

By adjunction we get:
\begin{equation}
\RHom ^{\bullet}(\pi ^{\ast}\Ee ,{\sf F}_{\ast}\Oo _X) =0.
\end{equation}

Lemma \ref{lem:Frobgeneratorforp>h} then implies that the
object $\pi ^{\ast}\Ee = 0$ in $\Dd ^{-}(X')$. Applying to it the
functor $\pi _{\ast}$ and using the projection formula, we get ${\rm
  R}^{\bullet}\pi _{\ast}\pi ^{\ast}\Ee=\Ee \otimes {\rm
  R}^{\bullet}\pi _{\ast}\Oo _{X'}=\Ee$, hence $\Ee = 0$, q.e.d.
\end{proof}

At the end of this section, let us give the simplest example. Consider
the projective space $\Pp ^n$. It is an old fact going back to
\cite{Harample} that the Frobenius pushforward of a line bundle on
$\Pp ^n$ decomposes into the direct sum of line bundles. The set of
line bundles constituting the decomposition of ${\sf F}_{\ast}\Oo
_{\Pp ^n}$ is easily seen to belong to $\Oo _{\Pp ^n},\Oo _{\Pp
  ^n}(-1),\dots,\Oo _{\Pp ^n}(-n)$. Using that cohomology of line
bundles on $\Pp ^n$ can at most be non-zero in one degree, we see that
$\Ext ^i({\sf F}_{\ast}\Oo _{\Pp ^n},{\sf F}_{\ast}\Oo _{\Pp
  ^n})=0$ for $i>0$. Further, if $p>n+1$ (the Coxeter number of the
root system of type ${\bf A}_n$), then Corollary
\ref{cor:Frobgeneratorforp>h} states that ${\sf F}_{\ast}\Oo _{\Pp
  ^n}$ is a generator in $\Dd ^{b}(\Pp ^n)$. It follows that every
line bundle $\Oo _{\Pp ^n}(-i)$ for $0\leq i\leq -n$ comes in the
decomposition of ${\sf F}_{\ast}\Oo _{\Pp ^n}$ with a non-zero
multiplicity. Indeed, if a line bundle was absent for some $i$, then ${\sf F}_{\ast}\Oo _{\Pp
  ^n}$ would not be a generator. Hence,  ${\sf F}_{\ast}\Oo _{\Pp
  ^n}$ is a tilting bundle on $\Pp ^n$ if $p>n+1$, and one recovers Beilinson's
full exceptional collection \cite{Beil}.


\section{A few lemmas}\label{subsec:genstat}


\subsection{$\Ext$-groups}\label{subsec:techlemmas}

Recall that for a variety $X$ the $n$-th Frobenius twist of $X$ is
denoted $X^{(n)}$. One has a morphism of $k$-schemes ${\sf F}^{n}:X\rightarrow X^{(n)}$.

\bigskip

Let $\pi : Y\rightarrow X^{(n)}$ be an arbitrary morphism. Consider the
cartesian square:


\begin{figure}[H]\label{fig:fig14}
$$\xymatrix @C5pc @R4pc {
 {\tilde Y}\ar[r]^{p_{2}} \ar@<-0.1ex>[d]^{p_{1}} & Y \ar@<-0.4ex>[d]^{\pi} \\
       X   \ar[r]^{{\sf F}^{n}} & X^{(n)}
 }
$$
\end{figure}

\begin{lemma}\label{lem:fibsquare1}
The fibered product ${\tilde Y}$ is isomorphic 
  to the left uppermost corner in the cartesian square:


\begin{figure}[H]
$$
\xymatrix @C5pc @R4pc {
 {\tilde Y}\ar[r] \ar@<-0.1ex>[d]^{i} & \Delta ^{(n)}
 \ar@<-0.4ex>[d]^{i_{\Delta ^{(n)}}} \\
       X\times Y   \ar[r]^{{\sf F}^{n}\times \pi} & X^{(n)}\times X^{(n)}
 }
$$
\end{figure}

where ${\Delta ^{(n)}}$ is the diagonal in $X^{(n)}\times X^{(n)}$. If 
  $\pi$ is flat, then one has an isomorphism of sheaves
  $i_{\ast}{\Oo _{\tilde Y}} = ({\sf F}^{n}\times
  \pi)^{\ast}({i _{\Delta ^{(n)}}}_{\ast}{\Oo _{\Delta ^{(n)}}})$.
\end{lemma}

\begin{proof} 
The isomorphism of two fibered products follows
from the definition of fibered product. The isomorphism of sheaves follows from flatness
  of the Frobenius morphism and from flat base change. 
\end{proof}

\begin{definition}\label{lem:frobneighb}
Let $Y^{(n)}\stackrel{i}{\hookrightarrow} X^{(n)}$ be a closed subscheme.
The fibered product $Y^{(n)}\times _{X^{(n)}} X$ as defined in the
diagram:

\begin{figure}[H]
$$
\xymatrix @C5pc @R4pc {
 Y' \times _{X^{(n)}} X\ar[r] \ar@<-0.1ex>[d] & Y^{(n)} \ar@<-0.4ex>[d]^{i} \\
 X\ar[r]^{{\sf F}^n} & X^{(n)}
 }
$$
\end{figure}

is called the $n$-th Frobenius neighbourhood of the subscheme
 $Y^{(n)}$ in $X$.
\end{definition}
Consider the cartesian square:

\begin{figure}[H]
$$
\xymatrix @C5pc @R4pc {
 {\tilde X}\ar[r]^{\pi _{2}} \ar@<-0.1ex>[d]^{\pi _{1}} & X
 \ar@<-0.4ex>[d]^{{\sf F}^n} \\
       X   \ar[r]^{{\sf F}^n} & X^{(n)}
 }
$$
\end{figure}

\begin{lemma}\label{lem:diagonal}
The fibered product ${\tilde X}$ is isomorphic to the $n$-th
  Frobenius neighbourhood of the diagonal $\Delta ^{(n)} \subset X^{(n)}\times
  X^{(n)}$.
\end{lemma}

\begin{proof}
Apply Lemma \ref{lem:fibsquare1} to $Y = X$ and $\pi = {\sf F}^n$.
\end{proof}

Recall that a right adjoint functor ${{\sf F}^n}^{!}(?)$ to ${\sf F}^n_{\ast}(?)$
is isomorphic to ${{\sf F}^n}^{\ast}(?)\otimes \omega _{X}^{1-p^n}$. We get:
\begin{eqnarray}\label{eq:ext-cohom}
& \Ext ^{k}_{X}({\sf F}^n_{\ast}\Oo _{X},{\sf F}^n_{\ast}\Oo _{X}) = \Ext
^{k}({\Oo _{X}},{{\sf F}^n}^{!}{\sf F}^n_{\ast}\Oo _{X}) = & \nonumber \\
& =  \Ext ^{k}({\Oo _{X}},{{\sf F}^n}^{\ast}{\sf F}^n_{\ast}{\Oo _{X}}\otimes \omega _{X}^{1-p^n}) =
{\rm H}^{k}(X,{{\sf F}^n}^{\ast}{\sf F}^n_{\ast}{\Oo _{X}}\otimes \omega _{X}^{1-p^n}).
\end{eqnarray}

\begin{lemma}\label{lem:isom}
There is an isomorphism of cohomology groups:
\begin{eqnarray}\label{eq:mainisom}
& {\rm H}^{k}(X,{{\sf F}^n}^{\ast}{\sf F}^n_{\ast}{\Oo _{X}}\otimes \omega _{X}^{1-p^{n}}) = & \nonumber\\ 
& = {\rm H}^{k}(X\times X,({{\sf F}^n\times {{\sf
    F}^n}})^{\ast}({i_{\Delta ^{(n)}}}_{\ast}{\Oo _{\Delta
    ^{(n)}}})\otimes ({\omega  _{X}^{1-p^{n}}}\boxtimes {\Oo _{X}})).
\end{eqnarray}
\end{lemma}

\begin{proof}
This lemma was proved in \cite{Sam1} (Lemma 2.3) for $n=1$. For
convenience of the reader, let us repeat the proof. Consider the
above cartesian square. By flat base change 
one gets an isomorphism of functors, the Frobenius morphism ${\sf F}_n$ being flat:
\begin{equation}\label{eq:basechangefrob}
{{\sf F}^n}^{\ast}{\sf F}^n_{\ast} = {\pi _{1}}_{\ast}{\pi _{2}}^{\ast}.
\end{equation}

Note that all the functors ${\sf F}^n_{\ast},{{\sf F}^n}^{\ast},{\pi
  _{1}}_{\ast},{\pi _{2}}^{\ast}$ are exact, the Frobenius morphism
  ${\sf F}_n$ being affine. The isomorphism (\ref{eq:basechangefrob}) implies an isomorphism of cohomology groups
\begin{equation}\label{eq:isom1}
{\rm H}^{k}(X,{{\sf F}^n}^{\ast}{\sf F}^n_{\ast}{\Oo _{X}}\otimes \omega _{X}^{1-p^{n}}) =
{\rm H}^{k}(X,{\pi _{1}}_{\ast}{\pi _{2}}^{\ast}{\Oo _{X}}\otimes \omega
_{X}^{1-p^{n}}).
\end{equation}

By projection formula the right-hand side group in (\ref{eq:isom1}) is isomorphic to 
${\rm H}^{k}({\tilde X},{\pi _{2}}^{\ast}{\Oo _{X}}\otimes {\pi _{1}}^{\ast}{\omega
  _{X}^{1-p^{n}}})$. Let $p_1$ and $p_2$ be the projections of $X\times X$ onto
  the first and the second component respectively, and let $\tilde i$
  be the embedding $\tilde X\hookrightarrow X\times X$. One sees that 
${\pi _1} = p_1\circ {\tilde i}, \ {\pi _2} = p_2\circ
{\tilde i}$. Hence an isomorphism of sheaves
\begin{equation}\label{eq:isom2}
{\pi
  _{2}}^{\ast}{\Oo _{X}}\otimes {\pi _{1}}^{\ast}{\omega
  _{X}^{1-p^{n}}} = {\tilde i}^{\ast}(p_{2}^{\ast}{\Oo _{X}}\otimes
  p_{1}^{\ast}{\omega _{X}^{1-p^{n}}}) = {\tilde i}^{\ast}({\omega _{X}^{1-p^{n}}}
  \boxtimes {\Oo _{X}}).
\end{equation}

From these isomorphisms and from the projection formula one
  gets 
\begin{eqnarray}\label{eq:isom3}
& {\rm H}^{k}({\tilde X},{\pi
  _{2}}^{\ast}{\Oo _{X}}\otimes {\pi _{1}}^{\ast}{\omega
  _{X}^{1-p^n}}) = {\rm H}^{k}({\tilde X},{\tilde i}^{\ast}({\omega _{X}^{1-p^n}}
  \boxtimes {\Oo _{X}})) = & \nonumber\\ 
& = {\rm H}^{k}(X\times X,{\tilde
  i}_{\ast}{\Oo _{\tilde X}}\otimes ({\omega _{X}^{1-p^n}}
  \boxtimes {\Oo _{X}})).
\end{eqnarray}

By Lemma \ref{lem:diagonal} the subscheme ${\tilde X}$
  is isomorphic to the $n$-th Frobenius neighbourhood of the diagonal ${\Delta ^{(n)}}$ in $X^n\times X^n$; thus
 
\begin{equation}\label{eq:isom4}
{\tilde i}_{\ast}{\Oo _{\tilde X}} 
= ({{\sf F}^n\times {\sf F}^n})^{\ast}({i_{\Delta ^{(n)}}}_{\ast}{\Oo _{\Delta ^{(n)}}}).
\end{equation}

Applying Lemma \ref{lem:fibsquare1} to $\pi = {\sf F}^n$ finishes the proof. 
\end{proof}

\begin{corollary}\label{cor:isomcor}
Let $\Ee _1$ and $\Ee _2$ be two vector bundles on $X$. There is an
isomorphism of cohomology groups:
\begin{eqnarray}
& \Ext ^{k}({\sf F}_{\ast}\Ee _1,{\sf F}_{\ast}\Ee _2) =
{\rm H}^{k}(X,{{\sf F}^n}^{!}{\sf F}^n_{\ast}(\Ee _2)\otimes \Ee _1^{\ast}) = & \nonumber\\
& = {\rm H}^{k}(X\times X,({\sf F}^n\times {\sf F}^n)^{\ast}({i_{\Delta
    ^{(n)}}}_{\ast}{\Oo _{\Delta ^{n}}})\otimes ((\Ee _1^{\ast}\otimes
\omega _{X}^{1-p^n})\boxtimes \Ee _2)) = & \nonumber\\
& =  {\rm H}^{k}(X\times X,({\sf F}^n\times {\sf F}^n)^{\ast}({i_{\Delta
    ^{(n)}}}_{\ast}{\Oo _{\Delta ^{n}}})\otimes (\Ee _2\boxtimes (\Ee _1^{\ast} \otimes \omega _{X}^{1-p^n})).
\end{eqnarray}

In particular, 
\begin{eqnarray}\label{eq:yetanotherisom}
& {\rm H}^{k}(X,{{\sf F}^n}^{\ast}{\sf F}^n_{\ast}{\Oo _{X}}\otimes \omega _{X}^{1-p^n}) = & \nonumber\\
& = {\rm H}^{k}(X\times X,({\sf F}^n\times {\sf F}^n)^{\ast}({i_{\Delta
    ^{(n)}}}_{\ast}{\Oo _{\Delta ^{n}}})\otimes ({\Oo  _{X}}\boxtimes
{\omega _{X}}^{1-p^n})).
\end{eqnarray}
\end{corollary}

\begin{proof}
Flat base change implies an isomorphism of functors
\begin{equation}\label{eq:theotherwayaround}
  {{\sf F}^n}^{\ast}{\sf F}^n_{\ast} = {\pi _{2}}_{\ast}{\pi _{1}}^{\ast}.
\end{equation}

Repeating verbatim the proof of Lemma \ref{lem:isom} one gets the
statement. 
\end{proof}

\subsection{$\Pp ^{1}$ - bundles}\label{subsec:P1bundles}

Assume given a smooth variety $S$ and a locally free sheaf $\Ee$ of
rank 2 on $S$. Let $X = \Pp _S(\Ee)$ be the projectivization of the bundle $\Ee$ 
and $\pi : X\rightarrow S$  the projection. Denote $\Oo _{\pi}(-1)$
the relative invertible sheaf such that ${\rm
  R}^{\bullet}\pi _{\ast}\Oo _{\pi}(1)=\Ee ^{\ast}$.
\begin{lemma}\label{lem:P^1Lemma}
For any $n\geq 1$ there is a short exact sequence of vector bundles on $X$:
\begin{equation}
0\rightarrow \pi ^{\ast}{\sf F}^n_{\ast}\Oo _{S} \rightarrow
{\sf F}^n_{\ast}\Oo _{X} \rightarrow \pi ^{\ast}({\sf F}^n_{\ast}({\sf
  D}^{p^n-2}\Ee \otimes \mbox{\rm det}\ \Ee)\otimes
\mbox{\rm det} \ \Ee ^{\ast})\otimes \Oo _{\pi}(-1)\rightarrow 0 .
\end{equation}
\end{lemma}

\begin{proof}
By Theorem \ref{th:Orlovth1},
the category $\Dd ^{b}(X)$ has a semiorthogonal decomposition:
\begin{equation}\label{eq:decomdercatP1bun}
{\rm D}^{b}(X) = \langle \Dd _{-1}, \Dd _{0} \rangle ,
\end{equation}

where $\Dd _{i}$ for $i=0,-1$ is a full subcategory of $\Dd ^{b}(X)$ that
consists of objects $\pi ^{\ast}({\mathcal F})\otimes \Oo _{\pi}(i)$,
for ${\mathcal F}\in \Dd ^{b}(S)$. Decomposition (\ref{eq:decomdercatP1bun})
means that for any object $A\in \Dd ^{b}(X)$ there is a distinguished triangle:
\begin{equation}\label{eq:exacttrianP1bun}
\dots \rightarrow \pi ^{\ast}{\rm R}^{\cdot}\pi _{\ast}A \rightarrow A \rightarrow
\pi ^{\ast}({\tilde A})\otimes \Oo _{\pi}(-1)\rightarrow \pi
^{\ast}{\rm R}^{\cdot}\pi _{\ast}A [1] \rightarrow \dots  .
\end{equation}

\smallskip

The object $\tilde A$ can be found by tensoring the triangle
(\ref{eq:exacttrianP1bun}) with $\Oo _{\pi}(-1)$ and applying the
functor ${\rm R}^{\bullet}\pi _{\ast}$ to the obtained triangle. Given that
${\rm R}^{\bullet}\pi _{\ast}\Oo _{\pi}(-1) = 0$, we get an isomorphism:
\begin{equation}\label{eq:isomforP1bun}
{\rm R}^{\bullet}\pi _{\ast}(A\otimes \Oo _{\pi}(-1)) \simeq
{\tilde A}\otimes {\rm R}^{\bullet}\pi _{\ast}\Oo _{\pi}(-2).
\end{equation}

One has ${\rm R}^{\bullet}\pi _{\ast}\Oo _{\pi}(-2) =
\mbox{det} \ \Ee[-1]$. Tensoring both sides of the isomorphism (\ref{eq:isomforP1bun}) with
$\mbox{det} \ \Ee ^{\ast}$, we get:
\begin{equation}\label{eq:tildeA}
{\tilde A} = {\rm R}^{\bullet}\pi _{\ast}(A\otimes \Oo _{\pi}(-1))\otimes \mbox{det} \
{\mathcal E}^{\ast}[1].
\end{equation}

Let now $A$ be the vector bundle ${\sf F}^n_{\ast}\Oo _X$. The triangle
(\ref{eq:exacttrianP1bun}) becomes in this case:
\begin{equation}\label{eq:trianforFrobdirimage}
\dots \rightarrow \pi ^{\ast}{\rm R}^{\bullet}\pi _{\ast}{\sf F}^n_{\ast}\Oo _{X} \rightarrow
{\sf F}^n_{\ast}\Oo _{X} \rightarrow
\pi ^{\ast}({\tilde A})\otimes \Oo _{\pi}(-1)\rightarrow \pi
^{\ast}{\rm R}^{\bullet}\pi _{\ast}{\sf F}^n_{\ast}\Oo _X [1] \rightarrow
\dots \quad .
\end{equation}

where ${\tilde A} = {\rm R}^{\bullet}\pi _{\ast}({\sf F}_{\ast}\Oo _{X}\otimes \Oo _{\pi}(-1))\otimes \mbox{det} \
{\mathcal E}^{\ast}[1]$. Recall that for a coherent sheaf $\Ff$ on
$X$ one has an isomorphism ${\rm R}^{i}\pi _{\ast}{\sf F}^n_{\ast}\Ff = {\sf F}^n_{\ast}{\rm
  R}^{i}\pi _{\ast}\Ff$ (see the remark after
Proposition \ref{prop:Frobprop}). Therefore,
\begin{equation}
{\rm R}^{\bullet}\pi _{\ast}{\sf F}^n_{\ast}\Oo _{X} =
{\sf F}^n_{\ast}{\rm R}^{\bullet}\pi _{\ast}\Oo _{X} = {\sf F}^n_{\ast}\Oo
_{S}.
\end{equation}

On the other hand, by the projection formula one has ${\rm R}^{\bullet}\pi _{\ast}({\sf
  F}^n_{\ast}\Oo _{X}\otimes \Oo _{\pi}(-1))$ = ${\rm R}^{\bullet}\pi
_{\ast}({\sf F}^n_{\ast}\Oo _{\pi}(-p^n))$ = ${\sf F}^n_{\ast}{\rm R}^{\bullet}\pi _{\ast}\Oo
_{\pi}(-p^n)$. The relative Serre duality for $\pi$ gives:
\begin{equation}
{\rm R}^{\bullet}\pi _{\ast}\Oo _{\pi}(-p^n) = {\sf D}^{p^n-2}\Ee\otimes \mbox{det} \ \Ee[-1].
\end{equation}

Let ${\tilde \Ee}$ be the vector bundle ${\sf D}^{p^n-2}\Ee\otimes \mbox{det} \ \Ee$.
Putting these isomorphisms together we see that the triangle
(\ref{eq:trianforFrobdirimage}) can be rewritten as follows:
\begin{equation}
\dots \rightarrow \pi ^{\ast}{\sf F}^n_{\ast}\Oo _{S} \rightarrow
{\sf F}^n_{\ast}\Oo _{X} \rightarrow \pi ^{\ast}({\sf F}^n_{\ast}{\tilde \Ee}\otimes
\mbox{det} \ \Ee ^{\ast})\otimes \Oo
_{\pi}(-1)\stackrel{[1]}\rightarrow \dots  .
\end{equation}

Therefore, the above distinguished triangle is in fact a
short exact sequence of vector bundles on $X$:
\begin{equation}\label{eq:frobshorexseqP1bun}
0\rightarrow \pi ^{\ast}{\sf F}^n_{\ast}\Oo _{S} \rightarrow
{\sf F}^n_{\ast}\Oo _{X} \rightarrow \pi ^{\ast}({\sf F}^n_{\ast}{\tilde \Ee}\otimes
\mbox{det} \ \Ee ^{\ast})\otimes \Oo _{\pi}(-1)\rightarrow 0.
\end{equation}

\end{proof}

\begin{remark} {\rm Using Corollary \ref{cor:isomcor} and Lemma
    \ref{lem:P^1Lemma}, in \cite{SamD-affflag} we prove the
    $\sf D$-affinity of the flag variety in type ${\bf B}_2$}.
\end{remark}

In the similar vein, let us prove yet another statement that we will use in the next section:

\begin{proposition}\label{prop:trivialmapS^2toD^2}
Let $\Ee$ be a rank two vector bundle over a smooth base $S$, and $\pi
:\Pp (\Ee)\rightarrow S$ the projection. Then one has the short exact
sequence on $S$:
\begin{equation}
0\rightarrow \mbox{\rm det} \ {\Ee^{\ast}}^{\otimes p^n}\rightarrow {{\sf F}^n}^{\ast}\Ee ^{\ast}\otimes
  {\sf S}^{p^n}\Ee ^{\ast}\rightarrow {\sf S}^{2p^n}\Ee
  ^{\ast}\rightarrow 0.
\end{equation}
\end{proposition}

\begin{proof}
Let $\Oo _{\pi}(1)$ be the relative invertible sheaf. Consider the relative Euler
  sequence:
\begin{equation}
0\rightarrow \Oo _{\pi}(-1)\rightarrow \pi ^{\ast}\Ee \rightarrow \Oo _{\pi}(1)\otimes \pi ^{\ast}\mbox{\rm det}
\ \Ee\rightarrow 0.
\end{equation}

Apply ${{\sf F}^n}^{\ast}$ to it:
\begin{equation}
0\rightarrow \Oo _{\pi}(-p^n)\rightarrow \pi ^{\ast}{{\sf F}^n}^{\ast}\Ee \rightarrow \Oo _{\pi}(p^n)\otimes \pi ^{\ast}\mbox{\rm det}
\ \Ee^{\otimes p^n}\rightarrow 0.
\end{equation}

Tensoring the above sequence with $\Oo _{\pi}(p^n)\otimes \pi
^{\ast}\mbox{\rm det} \ \Ee^{\otimes -p^n} $ and using the isomorphism
$\Ee \otimes \mbox{\rm det} \  \Ee = \Ee ^{\ast}$, we get:
\begin{equation}
0\rightarrow \pi ^{\ast}\mbox{\rm det} \ {\Ee^{\ast}}^{\otimes p^n}\rightarrow
\pi ^{\ast}{{\sf F}^n}^{\ast}\Ee ^{\ast}\otimes \Oo _{\pi}(p^n) \rightarrow \Oo _{\pi}(2p^n)\rightarrow 0.
\end{equation}

Finally, applying ${\rm R}^{\bullet}\pi _{\ast}$ and using the
projection formula, we get the statement.

\end{proof}

\begin{remark}
{\rm Tensoring the short exact sequence (\ref{eq:frobshorexseqP1bun}) with
  $\Oo _{\pi}(1)$ and applying the direct image ${\rm R}^{\cdot}\pi
  _{\ast}$ to this sequence, we get:
\begin{equation}
0\rightarrow {\sf F}_{\ast}\Oo _{S}\otimes {\rm R}^{\cdot}\pi _{\ast}\Oo _{\pi}(1) \rightarrow
 {\rm R}^{\cdot}\pi _{\ast}({\sf F}_{\ast}\Oo _{X}\otimes \Oo _{\pi}(1)) \rightarrow {\sf F}_{\ast}({\tilde \Ee})\otimes
\mbox{det} \ (\Ee ^{\ast})\rightarrow 0
\end{equation}

or, rather
\begin{equation}
0\rightarrow {\sf F}_{\ast}\Oo _{S}\otimes \Ee ^{\ast} \rightarrow
{\sf F}_{\ast}{\sf S}^{p}\Ee ^{\ast} \rightarrow {\sf
  F}_{\ast}({\tilde \Ee})\otimes \mbox{det} \ (\Ee ^{\ast})\rightarrow 0.
\end{equation}

Tensoring this sequence with $\mbox{det} \ (\Ee)$, we get:
\begin{equation}\label{eq:Frobpullbackofseqrank2}
0\rightarrow {\sf F}_{\ast}\Oo _{S}\otimes \Ee \rightarrow
{\sf F}_{\ast}{\sf S}^{p}\Ee \rightarrow {\sf F}_{\ast}{\tilde \Ee}\rightarrow 0.
\end{equation}

For a rank two vector bundle there is a well-known short
exact sequence (cf. \cite{Sam1}):
\begin{equation}\label{eq:seqrank2}
0\rightarrow {\sf F}^{\ast}\Ee \rightarrow {\sf S}^{p}\Ee \rightarrow {\tilde \Ee}\rightarrow 0,
\end{equation}

and we see that the sequence
(\ref{eq:Frobpullbackofseqrank2}) is obtained by applying the functor ${\sf F}_{\ast}$ to the sequence
(\ref{eq:seqrank2}).}

\end{remark}

\begin{remark}
{\rm Lemma \ref{lem:P^1Lemma} can be generalized for projective
  bundles of arbitrary rank. If $\Ee$ is a vector bundle of rank $n$
  over a scheme $S$ and $X = \Pp (\Ee)$ is the projective bundle,
  then there is a filtration on the bundle ${\sf F}_{\ast}\Oo _X$ with
  associated graded factors being of the form $\pi ^{\ast}\Ff
  _i\otimes \Oo _{\pi}(-i)$, where $\Ff _i$ are some vector bundles
  over $S$, and $0\leq i\leq n-1$. Let us work out an example of a
  vector bundle of rank 3. 
\begin{lemma}\label{lem:P2bundles}
Let $\Ee$ be a rank 3 vector bundle over a scheme $S$, and $X = {\Pp (\Ee)}$ the projective bundle. Then the bundle ${\sf F}^n_{\ast}\Oo
_X$ has a two-step filtration with associated graded factors being:
\begin{equation}
\pi ^{\ast}{\sf F}^n_{\ast}\Oo _S, \ \pi ^{\ast}{\Ee _1}\otimes \Oo _{\pi}(-1), \ \pi ^{\ast}({\sf
  F}^n_{\ast}{\Ee _2}\otimes \mbox{det} \ \Ee ^{\ast})\otimes \Oo _{\pi}(-2),
\end{equation}

where $\Ee _1$ is the cokernel of the embedding ${\sf F}^n_{\ast}\Oo
_S\otimes \Ee ^{\ast}\rightarrow {\sf F}^n_{\ast}{\sf S}^{p^n}\Ee ^{\ast}$
and $\Ee _2 = {\sf D}^{p^n-3}\Ee \otimes \mbox{\rm det}\Ee$.

\end{lemma}

\begin{proof}
Theorem \ref{th:Orlovth1} states that
  the category $\Dd ^{b}(X)$ has a semiorthogonal decomposition of
  three pieces:
\begin{equation}
\Dd ^{b}(X) = \langle \pi ^{\ast}\Dd ^{b}(S)\otimes \Oo
_{\pi}(-2),\pi ^{\ast}\Dd ^{b}(S)\otimes \Oo _{\pi}(-1),\pi ^{\ast}\Dd ^{b}(S)\rangle .
\end{equation}

This decomposition produces a distinguished triangle:
\begin{equation}\label{eq:1sttriangle}
\dots \rightarrow {\pi}^{\ast}{\sf F}^n_{\ast}\Oo _S\rightarrow {\sf
  F}^n_{\ast}\Oo _X\rightarrow {\mathcal A}\stackrel{[1]}\rightarrow \dots \quad ,
\end{equation}

where $\mathcal A$ is an object of $\Dd ^{b}(X)$ that, in turn, fits into a distinguished triangle:
\begin{equation}\label{eq:2ndtriangle}
\dots \pi ^{\ast}{\mathcal G}\otimes \Oo _{\pi}(-1)\rightarrow
{\mathcal A}\rightarrow \pi ^{\ast}\Ff \otimes \Oo
_{\pi}(-2)\stackrel{[1]}\rightarrow \dots \quad .
\end{equation}
\comment{
From the local description of the Frobenius morphism one sees that the
map ${\pi}^{\ast}{\sf F}^n_{\ast}\Oo _S\rightarrow {\sf F}^n_{\ast}\Oo
_X$ is an embedding of vector bundles.
}

Here $\mathcal G$ and $\Ff$ are objects of $\Dd ^{b}(X)$. To find
$\mathcal G$, tensor the triangle (\ref{eq:2ndtriangle}) with $\Oo _{\pi}(1)$ and
apply the functor ${\rm R}^{\bullet}\pi _{\ast}$. We get ${\mathcal
  G}={\rm R}^{\bullet}({\mathcal A}\otimes \Oo _{\pi}(1))$. The latter
can be found from the triangle (\ref{eq:1sttriangle}). Tensoring it
with $\Oo _{\pi}(1)$ and then applying ${\rm R}^{\bullet}\pi _{\ast}$,
and taking sheaf cohomology of the obtained triangle in $\Dd ^{b}(S)$,
we obtain a short exact sequence on $S$: 
\begin{equation}
0\rightarrow {\sf F}^n_{\ast}\Oo _S\otimes \Ee ^{\ast}\rightarrow {\sf
  F}^n_{\ast}{\sf S}^{p^n}\Ee ^{\ast}\rightarrow {\mathcal
  G}\rightarrow 0 ,
\end{equation}

where the first map is obtained by applying ${\sf F}^n_{\ast}$ to the
natural embedding of vector bundles ${{\sf F}^n}^{\ast}\Ee
^{\ast}\rightarrow {\sf S}^{p^n}\Ee ^{\ast}$. 

To find $\Ff$, tensor the triangle (\ref{eq:2ndtriangle}) with $\Oo _{\pi}(-1)$ and
apply the functor ${\rm R}^{\bullet}\pi _{\ast}$. We obtain an
isomorphism:
\begin{equation}
{\rm R}^{\bullet}\pi _{\ast}({\mathcal A}\otimes \Oo _{\pi}(-1)) = \Ff \otimes {\rm
  R}^{\bullet}\pi _{\ast}\Oo _{\pi}(-3).
\end{equation}

Now ${\rm R}^{\bullet}\pi _{\ast}\Oo _{\pi}(-3) = \mbox{det} \
  \Ee[-2]$. Tensoring the triangle (\ref{eq:1sttriangle}) with $\Oo
  _{\pi}(-1)$ and applying ${\rm R}^{\bullet}\pi _{\ast}$, we obtain an
  isomorphism:
\begin{equation}
{\rm R}^{\bullet}\pi _{\ast}({\mathcal A}\otimes \Oo _{\pi}(-1)) = {\rm
  R}^{\bullet}\pi _{\ast}({\sf F}^n_{\ast}\Oo _X\otimes \Oo _{\pi}(-1)) = {\rm
  R}^{\bullet}\pi _{\ast}({\sf F}^n_{\ast}\Oo _{\pi}(-p^n)) = {\sf F}^n_{\ast}{\rm
  R}^{\bullet}\pi _{\ast}\Oo _{\pi}(-p^n),
\end{equation}

and ${\rm R}^{\bullet}\pi _{\ast}\Oo _{\pi}(-p^n) = {\sf D}^{p^n-3}\Ee
\otimes \mbox{det} \Ee[-2]$ by the relative Serre duality. Hence,
\begin{equation}
\Ff = {\sf F}^n_{\ast}({\sf D}^{p^n-3}\Ee \otimes \mbox{det} \Ee)\otimes
\mbox{det} \Ee ^{\ast}.
\end{equation}
 This proves the statement.
\end{proof}
}
\end{remark}


\subsection{Blowups of surfaces}\label{subsec:Blowupsofsurfaces}

Let $X$ be a smooth variety and $Y$ its smooth subvariety of
codimension two. Consider the blow-up of $Y$ in $X$. Recall notations
from Subsection \ref{subsubsec:semdec}: there is a cartesian square:
\begin{figure}[H]
$$
\xymatrix @C5pc @R4pc {
 {\tilde Y}\ar[r]^{j} \ar@<-0.1ex>[d]^{p} & {\tilde X} \ar@<-0.4ex>[d]^{\pi} \\
       Y   \ar[r]^{i} & X
 }
$$
\end{figure}

\noindent Here $\tilde Y$ is the exceptional divisor. If ${\mathcal
  N}_{Y/X}$ is normal bundle to $Y$ in $X$, then the projection
$p$ is the projectivization of the bundle ${\mathcal N}_{Y/X}$. Let $\Oo
_{p}(-1)$ be the relative invertible sheaf.

\begin{lemma}\label{lem:frobblowupseq}
There is a short exact sequence:
\begin{equation}\label{eq:seqforblowup}
0\rightarrow \pi ^{\ast}{\sf F}_{\ast}\Oo _{X}\rightarrow {\sf F}_{\ast}\Oo
_{\tilde X}\rightarrow j_{\ast}(\Oo
_{p}(-1)\otimes p^{\ast}E)\rightarrow 0.
\end{equation}
\end{lemma}

\noindent Here $E$ is a coherent sheaf on $Y$ which fits into a short
exact sequence:
\begin{equation}\label{eq:459}
0\rightarrow E\rightarrow i^{\ast}\pi _{\ast}{\sf F}_{\ast}\Oo _{\tilde X}\rightarrow
{\rm R}^{0}p_{\ast}j^{\ast}{\sf F}_{\ast}\Oo _{\tilde X}\rightarrow 0.
\end{equation}

\begin{proof}
By Theorem \ref{th:Orlovth2}, the category 
$\Dd ^{b}(\tilde X)$ admits a semiorthogonal decomposition:
\begin{equation}
\Dd ^{b}(\tilde X) = \langle j_{\ast}(p^{\ast}\Dd ^{b}(X)\otimes \Oo
_{p}(-1)), \pi ^{\ast}\Dd ^{b}(X)\rangle .
\end{equation}

\noindent This means that there is a distinguished triangle:
\begin{equation}\label{eq:461}
\dots \rightarrow \pi ^{\ast}{\rm R}^{\bullet}\pi _{\ast}{\sf F}_{\ast}\Oo _{X}\rightarrow
{\sf F}_{\ast}\Oo _{\tilde X}\rightarrow j_{\ast}(\Oo _{p}(-1)\otimes
p^{\ast}E)\stackrel{[1]}\rightarrow \dots \quad .
\end{equation}

\noindent Consider the canonical morphism $\pi ^{\ast}{\rm R}^{\bullet}{\pi _{\ast}}{\sf F}_{\ast}\Oo
_{\tilde X}\rightarrow {\sf F}_{\ast}\Oo _{\tilde X}$. One has:
\begin{equation}
{\rm R}^{\bullet}\pi _{\ast}{\sf F}_{\ast}\Oo _{\tilde X} = {\sf F}_{\ast}{\rm R}^{\bullet}\pi _{\ast}\Oo
_{\tilde X} = {\sf F}_{\ast}\pi _{\ast}\Oo _{\tilde X} = {\sf F}_{\ast}\Oo _X. 
\end{equation}

\noindent Indeed, ${\rm R}^{\bullet}\pi _{\ast}\Oo _{\tilde X} = \Oo _X$. The morphism $\pi
^{\ast}{\sf F}_{\ast}\Oo _{X}\rightarrow {\sf F}_{\ast}\Oo _{\tilde X}$ is an
injective map of coherent sheaves at the generic point of $\tilde X$. 
Therefore it is an embedding of coherent sheaves, the sheaves $\pi
^{\ast}{\sf F}_{\ast}\Oo _X$ and ${\sf F}_{\ast}\Oo _{\tilde X}$ being locally free.\par

Taking sheaf cohomology ${\mathcal H}^{\bullet}$ of the sequence
(\ref{eq:461}) we see that the object $p^{\ast}\rm E$ has cohomology only
in degree zero, hence $\rm E$ is a coherent sheaf and the sequence
(\ref{eq:461}) in fact becomes a short exact sequence
(\ref{eq:seqforblowup}). To find the sheaf $\rm E$ let us apply the functor
$j^{\ast}$ to the sequence (\ref{eq:seqforblowup}):
\begin{equation}\label{eq:462}
 0\rightarrow {\sf L}^{1}j^{\ast}j_{\ast}(\Oo _{p}(-1)\otimes
p^{\ast}{\rm E})\rightarrow j^{\ast}\pi ^{\ast}{\sf F}_{\ast}\Oo _{X}\rightarrow
j^{\ast}{\sf F}_{\ast}\Oo _{\tilde X}\rightarrow \Oo _{p}(-1)\otimes
p^{\ast}{\rm E}\rightarrow 0.
\end{equation}

\noindent Recall that the normal bundle ${\mathcal N}_{{\tilde Y}/{\tilde
    X}}$ to $\tilde Y$ in $\tilde X$ is isomorphic to $\Oo _{\tilde
  Y}(\tilde Y) = \Oo _{p}(-1)$. Hence, the sheaf ${\sf L}^{1}j^{\ast}j_{\ast}(\Oo _{p}(-1)\otimes
p^{\ast}\rm E)$ is isomorphic to $\Oo _{p}(-1)\otimes p^{\ast}{\rm E}\otimes
\Oo _{\tilde Y}(-{\tilde Y}) = p^{\ast}\rm E$, and the sequence
(\ref{eq:462}) becomes
\begin{equation}\label{eq:463}
0\rightarrow p^{\ast}{\rm E}\rightarrow j^{\ast}\pi ^{\ast}{\sf F}_{\ast}\Oo _{X}\rightarrow
j^{\ast}{\sf F}_{\ast}\Oo _{\tilde X}\rightarrow \Oo _{p}(-1)\otimes
p^{\ast}{\rm E}\rightarrow 0.
\end{equation}

\noindent Applying now to the sequence (\ref{eq:463}) the functor
${\rm R}^{\bullet}p_{\ast}$ and taking into account that ${\rm R}^{\bullet}p_{\ast}\Oo _{p}(-1)=0$ we
get a short exact sequence 
\begin{equation}
0\rightarrow {\rm E}\rightarrow {\rm R}^{0}p_{\ast}j^{\ast}\pi ^{\ast}{\sf F}_{\ast}\Oo
_X\rightarrow {\rm R}^{0}p_{\ast}j^{\ast}{\sf F}_{\ast}\Oo _{\tilde X}\rightarrow 0.
\end{equation}

\noindent Now ${\rm R}^{0}p_{\ast}j^{\ast}\pi ^{\ast}{\sf F}_{\ast}\Oo _{X} =
{\rm R}^{0}p_{\ast}p^{\ast}i^{\ast}{\sf F}_{\ast}\Oo _X = i^{\ast}{\sf F}_{\ast}\Oo _X =
i^{\ast}\pi _{\ast}{\sf F}_{\ast}\Oo _{\tilde X}$, and the sequence
(\ref{eq:459}) follows.
\end{proof}

Consider a particular case when $X$ is a smooth surface and $Y$ is a
point $y\in X$. Let $\tilde X$ be the blown-up surface and $l$ be the
exceptional divisor, $l = \Pp ^{1}$. 
\begin{corollary}\label{cor:blowupforsurf}
There is a short exact sequence:
\begin{equation}\label{eq:465}
0\rightarrow \pi ^{\ast}{\sf F}_{\ast}\Oo _{X}\rightarrow {\sf F}_{\ast}\Oo
_{\tilde X}\rightarrow j_{\ast}\Oo _{l}(-1)^{\oplus \frac{p(p-1)}{2}}\rightarrow 0.
\end{equation}
\end{corollary}

\begin{proof}
The category $\Dd ^{b}(y)$ is equivalent to $\Dd ^{b}(\mbox{Vect}-k)$,
since $y$ is a point. Hence, we just need to compute the multiplicity
of the sheaf $j_{\ast}\Oo _l(-1)$ in the sequence (\ref{eq:465}) or
the rank of vector space $E$. 
This multiplicity is equal to the corank of the morphism of sheaves
$\pi ^{\ast}{\sf F}_{\ast}\Oo _X\rightarrow {\sf F}_{\ast}\Oo _{\tilde X}$ at the point $y$.

\begin{proposition}
The corank is equal to $\frac{p(p-1)}{2}$.
\end{proposition}

\begin{proof}
Choose the local coordinates $x, y$ on $\tilde X$. Then $x, xy$ are
the local coordinates on $X$. The stalk of the sheaf $\pi
^{\ast}{\sf F}_{\ast}\Oo _X$ at $y$ is then $k[x,y]/(x^{p},(xy)^{p})$
whereas the stalk of the sheaf ${\sf F}_{\ast}\Oo _{\tilde X}$ at $y$ is
$k[x,y]/(x^{p},y^{p})$. We see now that the cokernel of the map
$k[x,y]/(x^{p},(xy)^{p})\rightarrow k[x,y]/(x^{p},y^{p})$ consists of
monomials $x^{a}y^{b}$ such that $0\leq a < b < p$, hence the statement.
\end{proof}

Corollary \ref{cor:blowupforsurf} is proven.
\end{proof}


\subsection{Cohomology of the Frobenius neighborhoods}\label{subsec:cohofFrobneigh}


Let $X$ be a smooth variety over $k$. To compute the 
    cohomology groups ${\rm H}^{i}(X',{\mathcal End}({{\sf
    F}}_{\ast}\Oo _X))$ we will use the properties of sheaves
    ${\rm D}_X$ from Section \ref{subsec:Berthdiffoper}. Keeping
    the previous notation, we get:
\begin{eqnarray}\label{eq:chainofisom}
& {\rm H}^{j}(X',{\mathcal End}({\sf F}_{\ast}\Oo _X)) = 
{\rm H}^{j}(X',i^{\ast}\DD _X) = 
{\rm H}^{j}({\rm T}^{\ast}(X'),i_{\ast}i^{\ast}\DD _X) = & \nonumber \\
& = {\rm H}^{j}({\rm T}^{\ast}(X'),\DD _X\otimes i_{\ast}\Oo _{X^{'}}), 
\end{eqnarray}

the last isomorphism in (\ref{eq:chainofisom}) follows from
  the projection formula. Recall the projection $\pi \colon {\rm
    T}^{\ast}(X^{'})\rightarrow X^{'}$. Consider the bundle $\pi ^{\ast}\T
  ^{\ast}_{X^{'}}$. There is a tautological section $s$ of this bundle
such that the zero locus of $s$ coincides with $X^{'}$. Hence, one obtains the Koszul resolution:
\begin{equation}\label{eq:Kosres}
0\rightarrow \dots \rightarrow \bigwedge^{k}(\pi ^{\ast}\T _{X^{'}}) \rightarrow
\bigwedge^{k-1}(\pi ^{\ast}\T _{X^{'}})\rightarrow \dots \rightarrow \Oo
_{{\rm T}^{\ast}(X^{'})}\rightarrow i_{\ast}\Oo _{X^{'}}\rightarrow 0.
\end{equation}

Tensor the resolution (\ref{eq:Kosres}) with the sheaf
$\DD _X$. The rightmost cohomology group in (\ref{eq:chainofisom})
can be computed via the above Koszul resolution. 
\begin{lemma}\label{lem:cohofFrobneigh}
\ \ Fix \ \ $k\geq 0$. \ \ For \ any \ $j\geq 0$, \ if \ ${\rm H}^{j}({\rm
  T}^{\ast}(X),{\sf F}^{\ast}\bigwedge ^{k}(\pi ^{\ast}\T _{X'})) = 0$,
  \ then \ ${\rm H}^{j}({\rm T}^{\ast}(X'),\DD _X\otimes \bigwedge ^{k}(\pi ^{\ast}\T _{X'})) = 0$. 
\end{lemma}

\begin{proof}
Denote $C^{k}$ the sheaf $\DD _X \otimes \bigwedge ^{k}(\pi ^{\ast}\T _{X'})$.
Take the direct image of $C^k$ with respect to $\pi$. Using the
projection formula we get:
\begin{eqnarray}
& {\rm R}^{\bullet}\pi _{\ast}C^{k} = {\rm
  R}^{\bullet}\pi _{\ast}(\DD _X\otimes \bigwedge ^{k}(\pi ^{\ast}\T
  _{X'})) = \pi _{\ast}(\DD _X\otimes \bigwedge ^{k}(\pi ^{\ast}\T
  _{X'})) & \nonumber \\
& = {\sf F}_{\ast}{\rm D} _X\otimes \bigwedge ^{k}(\T _{X'}), 
\end{eqnarray}

the morphism $\pi$ being affine. Hence,
\begin{equation}
{\rm H}^{j}({\rm T}^{\ast}(X'),\DD _X\otimes
\bigwedge ^{k}(\pi ^{\ast}\T _{X'})) = {\rm H}^{j}(X', {\sf
 F}_{\ast}{\rm D} _X\otimes \bigwedge ^{k}(\T _{X'})).
\end{equation}

The sheaf ${\sf F}_{\ast}{\rm D} _X\otimes
\bigwedge ^{k}(\T _{X'})$ is equipped with a filtration that is induced
  by the filtration on ${\sf F}_{\ast}{\rm D}_X$, the associated sheaf
  being isomorphic to ${\rm gr}({\sf F}_{\ast}{\rm D} _X)\otimes \bigwedge ^{k}(\T
_{X'}) = {\sf F}_{\ast}\pi _{\ast}\Oo _{{\rm T}^{\ast}(X)}\otimes
\bigwedge ^{k}(\T _{X'})$. Clearly, for $j\geq 0$ 
\begin{equation}
{\rm H}^{j}(X',{\rm gr}({\sf F}_{\ast}{\rm D} _X)\otimes
\bigwedge ^{k}(\T _{X'})) = 0 \quad \Rightarrow \quad {\rm H}^{j}(X',{\sf F}_{\ast}{\rm D}_{X}\otimes
\bigwedge ^{k}(\T _{X'})) = 0. 
\end{equation}

There are isomorphisms:
\begin{eqnarray}
& {\rm H}^{j}(X',{\sf F}_{\ast}\pi _{\ast}\Oo _{{\rm T}^{\ast}(X)}\otimes
\bigwedge ^{k}(\T _{X'})) = {\rm H}^{j}(X,\pi _{\ast}\Oo _{{\rm
  T}^{\ast}(X)}\otimes {\sf F}^{\ast}\bigwedge ^{k}(\T _{X'})) = & \nonumber \\
& = {\rm H}^{j}({\rm T}^{\ast}(X),\pi ^{\ast}{\sf F}^{\ast}\bigwedge ^{k}(\T
_{X'})),
\end{eqnarray}

the \ last \ isomorphism \ following  \ from \ the \ projection \
formula. \ Finally, \ ${\rm H}^{j}({\rm T}^{\ast}(X),\pi ^{\ast}{\sf F}^{\ast}\bigwedge ^{k}(\T
_{X'})) = {\rm H}^{j}({\rm T}^{\ast}(X),{\sf F}^{\ast}\bigwedge ^{k}(\pi
^{\ast}\T _{X'}))$, hence the statement of the lemma. 
\end{proof}

\begin{remark}
Assume that for a given $X$ one has ${\rm H}^{j}({\rm T}^{\ast}(X),{\sf F}^{\ast}\bigwedge
^{k}(\pi ^{\ast}\T _{X'})) = 0$ for $j>k\geq 0$. The \ spectral \ sequence \
${\rm E}_1^{p,q} = {\rm H}^{p}({\rm T}^{\ast}(X'),\DD _X\otimes
\bigwedge ^q(\pi ^{\ast}\T _{X'}))$ \ converges \ to \ ${\rm H}^{p-q}({\rm
  T}^{\ast}(X'),\DD _X\otimes i_{\ast}\Oo _{X^{'}})$.  Lemma \ref{lem:cohofFrobneigh} and the resolution
(\ref{eq:Kosres}) then imply:

\begin{equation}
{\rm H}^{j}({\rm
  T}^{\ast}(X'),\DD _X\otimes i_{\ast}\Oo _{X^{'}}) = {\rm H}^{j}(X',{\mathcal End}({\sf F}_{\ast}\Oo _X)) = 0
\end{equation}

for $j>0$, and 

\begin{equation}
{\rm H}^{j}({\rm  T}^{\ast}(X),{\sf F}^{\ast}i_{\ast}\Oo _{X'}) = 0
\end{equation}
 for $j>0$.
\end{remark}


\section{Flag varieties}\label{sec:flagvar}


In this section we show some applications of the above results. 

\subsection{Flag variety in type ${\bf A}_2$}

Consider the group ${\bf SL}_3$ over $k$ and
the flag variety ${\bf SL}_3/{\bf B}$. In \cite{Haas} it was proved
(Theorem 4.5.4 in {\it loc.cit.})
that the sheaf of differential operators on ${\bf SL}_3/{\bf B}$ has
vanishing higher cohomology (recall that for flag varieties this
implies the $\sf D$-affinity, see (\cite{Haas})). Below we give a different proof of this
vanishing theorem.

\begin{theorem}\label{th:SL_3}
Let $X$ be the flag variety ${\bf SL}_3/{\bf B}$.
Then $\Ext ^{i}({{\sf F}_{m}}_{\ast}\Oo _X,{{\sf F}_{m}}_{\ast}\Oo
_X)=0$ for $i>0$ and $m\geq 1$.
\end{theorem}

\begin{corollary}
${\rm H}^{i}(X,\D _X)=0$ for $i>0$.
\end{corollary}

\begin{proof}
Recall (see Subsection \ref{subsec:diffoper}) that the sheaf $\D _X$ is the direct limit of
sheaves of matrix algebras:
\begin{equation}
\D _X = \bigcup _{n\geq 1} {\mathcal End}_{\Oo _X}({{\sf
    F}_n}_{\ast}\Oo _X).
\end{equation}

Clearly, for some $i$ and $n\geq 1$ the vanishing 
$\Ext ^{i}({{\sf F}_{n}}_{\ast}\Oo _X,{{\sf F}_{n}}_{\ast}\Oo _X) = {\rm H}^{i}(X,{\mathcal End}({{\sf
    F}_n}_{\ast}\Oo _X)) = 0$ implies ${\rm H}^{i}(X,\D _X) = 0$.
\end{proof}

\begin{proof}
The flag variety ${\bf SL}_3/{\bf B}$ is isomorphic to an incidence
variety. For convenience of the reader, recall some facts about
incidence varieties. Let $V$ be a vector space of dimension $n$. The
incidence variety $X_n$ is the set of pairs $X_n\colon = (l\subset {\rm H}\subset V)$,
where $l$ and $\rm H$ are a line and a hyperplane in $V$,
respectively. The variety $X_n$ is fibered over $\Pp (V)$ and $\Pp (V^{\ast})$:

\begin{figure}[H]\label{fig:fig18}
$$
\xymatrix @C5pc @R4pc {
 & X_n \ar[dl]_{p} \ar@<-0.1ex>[dr]^{\pi} &   \\
 \Pp (V)           &                    &         \Pp (V^{\ast})
 }
$$
\end{figure}

Let $0\subset {\mathcal U}_1\subset {\mathcal U}_{n-1}\subset V\otimes \Oo _X$
be the tautological flag on $X_n$. The projection $\pi$ is 
projectivization of the bundle $\Omega ^{1}(1)$ on $\Pp
(V^{\ast})$. Let $\Oo _{p}(-1)$ and $\Oo
_{\pi}(-1)$ be the relative tautological line bundles with respect to
projections $p$ and $\pi$, respectively. Note that ${\mathcal U}_1 =
p^{\ast}\Oo (-1) = \Oo _{\pi}(-1), {\mathcal U}_{n-1} = \pi
^{\ast}\Omega ^{1}(1)$. Let ${\pi}/l$ be the quotient bundle: 
\begin{equation}\label{eq:tautseq}
0\rightarrow p^{\ast}\Oo (-1)\rightarrow \pi ^{\ast}\Omega
^{1}(1)\rightarrow {\pi}/l\rightarrow 0.
\end{equation}

Denote $\Oo (i,j)$ the line bundle $p^{\ast}\Oo
(i)\otimes \pi ^{\ast}\Oo (j)$. The canonical line bundle $\omega
_{X}$ is isomorphic to $\Oo (-n,-n)$. To compute the $\Ext$-groups,
let us apply Lemma \ref{cor:isomcor}. Recall that this lemma states an
isomorphism of the following groups:
\begin{equation}
\Ext ^{i}_X({{\sf F}_m}_{\ast}\Oo _X,{{\sf F}_m}_{\ast}\Oo _X) 
= {\rm H}^{k}(X\times X,({\sf F}_m\times {\sf F}_m)^{\ast}({i_{\Delta
    ^{(m)}}}_{\ast}{\Oo _{\Delta ^{m}}})\otimes ({\Oo  _{X}}\boxtimes
{\omega _{X}}^{1-p^m})).
\end{equation}

For incidence varieties, however, there is a nice resolution
of the sheaf $i_{\ast}\Oo _{\Delta}$, the Koszul resolution
(\cite{Kap}, Proposition 4.17). Recall briefly its construction.
Consider the following double complex of
sheaves $C^{\bullet,\bullet}$ on $X_n\times X_n$:

\begin{figure}[H]\label{fig:fig19}
$$
\xymatrix @C5pc @R4pc {
 \dots\ar[r] & \Psi _{1,0}\boxtimes \Oo (-1,0)\ar[r] & \Oo _{X_n\times X_n} &  \\
 \dots\ar[r] & \Psi _{1,1}\boxtimes \Oo (-1,-1)\ar[r]\ar[u] & \Psi _{0,1}\boxtimes \Oo (0,-1)\ar[u] & \\
 & \vdots \ar[u] &  \vdots \ar[u]  &
 }
$$
\end{figure}

The total complex of $C^{\bullet,\bullet}$ is a left
resolution of the structure sheaf of the diagonal $\Delta \subset
X_n\times X_n$. Truncate $C^{\bullet,\bullet}$, deleting all terms
except those belonging to the intersection of the first $n$ rows (from
0-th up to $(n-1)$-th) and the first $n-1$ columns, and consider the
convolution of the remaining double complex. Denote $\tilde C^{\bullet}$
the convolution. The truncated complex has
only two non-zero cohomology: ${\mathcal H}^{0} = \Oo _{\Delta}$, and
${\mathcal H}^{-2(n-1)}$. The latter cohomology can be explicitly
described:

\begin{equation}\label{eq:kerresofdiagincquad}
{\mathcal H}^{-2(n-1)} = \bigoplus _{i=0}^{n-1}\wedge ^{i}({\pi}/l)(-1,0)\boxtimes \wedge ^{i}({\pi}/l)^{\ast}(-n+1,-n).
\end{equation}

It follows that there is
the following distinguished triangle ($\sigma _{\geq}$ stands for the stupid truncation): 
\begin{equation}\label{disttrianforpartialflags}
\dots \rightarrow {\mathcal H}^{-2n+2}[2n-2]\rightarrow \sigma _{\geq
-2n+2}({\tilde C}^{\bullet})\rightarrow i_{\ast}\Oo
_{\Delta}\rightarrow {\mathcal H}^{-2n+2}[2n-1]\rightarrow \dots \quad .
\end{equation}

Let us come back to the case of ${\bf SL}_3/{\bf B} = X_3 = X$.
Using the above triangle and Lemma \ref{cor:isomcor}, we can
prove Theorem \ref{th:SL_3} almost immediately. In this case, the
above triangle looks as follows:
\begin{equation}\label{disttrianforpartialflags}
\dots \rightarrow {\mathcal H}^{-2}[2]\rightarrow \sigma _{\geq
-2}({\tilde C}^{\bullet})\rightarrow i_{\ast}\Oo
_{\Delta}\rightarrow {\mathcal H}^{-2}[3]\rightarrow \dots  ,
\end{equation}

where the truncated complex $\sigma _{\geq
  -2}({\tilde C}^{\bullet})$ is quasiisomorphic to:
\begin{equation}\label{eq:SL3truncated}
0\rightarrow \Psi _{1,1}\boxtimes \Oo (-1,-1)\rightarrow 
\Psi _{1,0}\boxtimes \Oo (-1,0)\oplus \Psi _{0,1}\boxtimes \Oo
(0,-1)\rightarrow \Oo _{X\times X}\rightarrow 0,
\end{equation}

and there is an isomorphism

\begin{equation}\label{eq:kerresofdiagincquad}
{\mathcal H}^{-2} = \Oo (-1,0)\boxtimes \Oo (-1,-2)\oplus
{\pi/l}\otimes \Oo (-1,0)\boxtimes ({\pi}/l)^{\ast}\otimes \Oo (-1,-2).
\end{equation}

We need to compute the groups ${\rm H}^{k}(X\times X,({\sf F}_m\times {\sf F}_m)^{\ast}({i_{\Delta
    ^{(m)}}}_{\ast}{\Oo _{\Delta ^{m}}})\otimes ({\Oo  _{X}}\boxtimes
{\omega _{X}}^{1-p^m}))$. Apply the functor $({\sf F}_m\times {\sf
  F}_m)^{\ast}$ to the triangle (\ref{disttrianforpartialflags}) and
tensor it then with the sheaf $({\Oo  _{X}}\boxtimes
{\omega _{X}}^{1-p^m}))$. Let us first prove that ${\mathbb
  H}^{i}(X\times X,({\sf F}_m\times {\sf F}_m)^{\ast}(\sigma _{\geq
-2}({\tilde C}^{\bullet}))\boxtimes ({\Oo  _{X}}\boxtimes
{\omega _{X}}^{1-p^m}))$ = 0 for $i>0$. Recall that $\omega _X = \Oo
(-2,-2)$. The sheaves $\Psi _{i,j}$ have right resolutions
consisting of direct sums of ample line bundles and the sheaf $\Oo _X$
(\cite{Kap}). This implies 
that ${\rm H}^{i}(X,{\sf F}_m^{\ast}(\Psi _{1,0})) = {\rm H}^{i}(X,{\sf F}_m^{\ast}(\Psi _{0,1})) = 0$
for $i>1$ and ${\rm H}^{i}(X,{\sf F}_m^{\ast}(\Psi _{1,1})) = 0$ for $i>2$. Along the second argument in
(\ref{eq:SL3truncated}) we get, after tensoring it with $\omega _{X}^{1-p^m}$, ample
line bundles. Ample line bundles have no higher cohomology by the
Kempf theorem (\cite{Ke}). The spectral sequence then gives ${\mathbb
  H}^{i}(X\times X,({\sf F}_m\times {\sf F}_m)^{\ast}(\sigma _{\geq
-2}({\tilde C}^{\bullet}))\boxtimes ({\Oo  _{X}}\boxtimes
{\omega _{X}}^{1-p^m}))$ = 0 for $i>0$.\par

Note that ${\pi}/l$ is a line bundle. One sees that ${\pi}/l = \Oo
(1,-1)$, hence the sheaf ${\mathcal H}^{-2}$ is isomorphic to 
$\Oo (-1,0)\boxtimes \Oo (-1,-2)\oplus \Oo (0,-1)\boxtimes \Oo
(-2,-1)$. Thus
\begin{eqnarray}
& ({\sf F}_m\times {\sf F}_m)^{\ast}{\mathcal H}^{-2}\otimes ({\Oo
 _{X}}\boxtimes {\omega _{X}}^{1-p^m}) = \\
& = \Oo (-p^m,0)\boxtimes \Oo (-p^m,-2p^m)\otimes \omega _X^{1-p^m}\oplus \Oo
(0,-p^m)\boxtimes \Oo (-2p^m,-p^m)\otimes \omega _X^{1-p^m}. &\nonumber
\end{eqnarray}
Let us prove that the latter bundle has only one non-vanishing cohomology, namely
${\rm H}^{3}$. The Serre duality gives:
\begin{equation}\label{eq:partialflagisom1}
{\rm H}^{i}(X,{\sf F}_m^{\ast}\Oo (-1,0)) = {\rm H}^{3-i}(X,{\sf F}_m^{\ast}\Oo (-1,-2)\otimes
\omega _X^{1-p^m}),
\end{equation}

and 
\begin{equation}\label{eq:partialflagisom2}
{\rm H}^{i}(X,{\sf F}_m^{\ast}\Oo (0,-1)) = {\rm H}^{3-i}(X,{\sf F}_m^{\ast}\Oo (-2,-1)\otimes
\omega _X^{1-p^m}).
\end{equation}

It is therefore sufficient to show that the left-hand sides
in both (\ref{eq:partialflagisom1}) and (\ref{eq:partialflagisom2})
are non-zero only for one value of $i$. Consider for example the cohomology group
${\rm H}^{i}(X,{\sf F}_m^{\ast}\Oo (-1,0))$. The line bundle $\Oo (-1,0)$ is
isomorphic to $p^{\ast}\Oo (-1)$, hence ${\rm H}^{i}(X,{\sf F}_m^{\ast}\Oo
(-1,0)) = {\rm H}^{i}(\Pp ^2,{\sf F}_m^{\ast}\Oo (-1)) = {\rm H}^{i}(\Pp
^2,\Oo (-p^m))$. The line bundle $\Oo (-p^m)$ on $\Pp ^2$ is either
acyclic (for $p=2$ and $m=1$) or has only one non-zero cohomology group in top
degree. The K\"unneth formula now gives that 
${\rm H}^{i}(X\times X,({\sf F}_m\times {\sf F}_m)^{\ast}{\mathcal H}^{-2}\otimes ({\Oo
_{X}}\boxtimes {\omega _{X}}^{1-p^m})) = 0$ for $i\neq 3$. Remembering
the distinguished triangle (\ref{disttrianforpartialflags}), we get the proof.
\end{proof}

\subsection{Flag variety in type ${\bf B}_2$}

Let $\sf V$ be a symplectic vector space of dimension 4 over $k$. Let $\bf G$ be the
symplectic group ${\bf Sp}_4$ over $k$; the root system of $\bf G$ is of type ${\bf
  B}_2$. Let $\bf B$ be a Borel subgroup of $\bf G$. Consider the flag variety ${\bf G}/{\bf B}$. The group $\bf G$ has two parabolic subgroups ${\bf P}_{\alpha}$ and ${\bf
  P}_{\beta}$ that correspond to the simple roots $\alpha$ and
$\beta$, the root $\beta$ being the long root. 
The homogeneous spaces ${\bf G}/{{\bf P}_{\alpha}}$ and  ${\bf G}/{{\bf
    P}_{\beta}}$ are isomorphic to the 3-dimensional quadric ${\sf Q}_3$ and $\Pp
^3$, respectively. Denote $q$ and $\pi$ the two projections of ${\bf G}/{\bf B}$ onto ${\sf Q}_3$ and $\Pp
^3$. The line bundles on ${\bf G}/{\bf B}$ that correspond to the fundamental weights $\omega _{\alpha}$ and
$\omega _{\beta}$ are isomorphic to $\pi ^{\ast}\Oo _{\Pp ^3}(1)$ and
$q^{\ast}\Oo _{{\sf Q}_3}(1)$, respectively. The canonical line bundle
$\omega _{{\bf G}/{\bf B}}$ corresponds to the weight $-2\rho =
-2(\omega _{\alpha} + \omega _{\beta})$ and is isomorphic to
$\pi ^{\ast}\Oo _{\Pp ^3}(-2)\otimes q^{\ast}\Oo _{{\sf Q}_3}(-2)$.
The projection $\pi$ is the projective
bundle over $\Pp ^3$ associated to a rank two vector bundle $\sf N$
over $\Pp ^3=\Pp ({\sf V})$, and the projection $q$ is the projective bundle associated to the spinor bundle $\Uu
_2$ on ${\sf Q}_3$. The bundle $\sf N$ is symplectic, that is there is
a non-degenerate skew-symmetric pairing $\wedge ^2{\sf N}\rightarrow
\Oo _{\Pp ^3}$ that is induced by the given symplectic structure on
$\sf V$. There is a short exact sequence on $\Pp ^3$:
\begin{equation}
0\rightarrow \Oo _{\Pp ^3}(-1)\rightarrow \Omega _{\Pp ^3}^1(1)\rightarrow {\sf
  N}\rightarrow 0,
\end{equation}

while the spinor bundle $\Uu _2$, which is also isomorphic to the restriction 
of the rank two universal bundle on ${\rm Gr}_{2,4}={\sf Q}_4$ to
${\sf Q}_3$, fits into a short exact sequence on ${\sf Q}_3$:
\begin{equation}\label{seq:universalseqonQ_3}
0\rightarrow \Uu _2\rightarrow {\sf V}\otimes \Oo _{{\sf
    Q}_3}\rightarrow \Uu _2^{\ast}\rightarrow 0.
\end{equation}

\begin{theorem}\label{th:B_2}
Let the characteristic of $k$ be odd. Then $\Ext ^{i}({\sf F}_{\ast}\Oo _{{\bf G}/{\bf B}},{\sf F}_{\ast}\Oo _{{\bf G}/{\bf B}})=0$ for $i>0$.
\end{theorem}

\begin{proof}
Let $\Oo _{\pi}(1)$ and $\Oo _{q}(1)$ be the relative line bundles
for the projections $\pi$ and $q$, respectively. In particular,
$q_{\ast}\Oo _{q}(1)=\Uu _2^{\ast}$. Note that $\Oo _{q}(1) = \pi
  ^{\ast}\Oo _{\Pp ^3}(1)$. Denote $\mathcal T$ the tangent bundle of ${\bf
  G}/{\bf B}$. By Lemma \ref{lem:cohofFrobneigh}, the statement will follow if we show
that the cohomology groups 
\begin{equation}\label{eq:FrNzerosection}
{\rm H}^{i}({\bf G}/{\bf B},\wedge ^{k}{\sf F}^{\ast}{\mathcal T}\otimes {\sf
  S}^{\bullet}{\mathcal T}) = {\rm H}^{i}({\bf G}/{\bf B},\wedge ^{k}{\mathcal
  T}\otimes {\sf F}_{\ast}{\sf S}^{\bullet}{\mathcal T})
\end{equation}

are zero for $i>k$. Consider the short exact sequence:
\begin{equation}\label{eq:tangseqforproj_pi}
0\rightarrow {\mathcal T}_{\pi}\rightarrow {\mathcal T}\rightarrow \pi
^{\ast}{\mathcal T}_{\Pp ^3}\rightarrow 0.
\end{equation}

Here ${\mathcal T}_{\pi}=\Oo _{\pi}(2)$ is the relative tangent sheaf with respect to the
  projection $\pi$. It is a line bundle on ${\bf G}/{\bf B}$ that
corresponds to the long root $\beta$. One has the Euler sequence on
  $\Pp ^3$:
\begin{equation}
0\rightarrow \Oo _{{\bf G}/{\bf B}}\rightarrow {\sf V}\otimes \pi
^{\ast}\Oo _{\Pp ^3}(1)\rightarrow \pi ^{\ast}{\mathcal T}_{\Pp ^3}\rightarrow 0.
\end{equation}

Note that there is an isomorphism of line bundles: $\Oo _{\pi}(2) = q^{\ast}\Oo (2)\otimes
\pi ^{\ast}\Oo _{\Pp ^3}(-2)$. Consider the short exact sequences for exterior powers of $\mathcal T$
that are obtained from the sequence (\ref{eq:tangseqforproj_pi}):
\begin{equation}\label{eq:extsqoftangentseq}
0\rightarrow \pi ^{\ast}{\mathcal T}_{\Pp ^3}\otimes \Oo
_{\pi}(2)\rightarrow \wedge ^2{\mathcal T}\rightarrow \wedge ^2\pi ^{\ast}{\mathcal T}_{\Pp ^3}\rightarrow 0,
\end{equation}

and

\begin{equation}
0\rightarrow \wedge ^2\pi ^{\ast}{\mathcal T}_{\Pp ^3}\otimes \Oo
_{\pi}(2)\rightarrow \wedge ^3{\mathcal T}\rightarrow \wedge ^3\pi ^{\ast}{\mathcal T}_{\Pp ^3}\rightarrow 0.
\end{equation}

The rest of the proof is broken up into a series
of propositions. We start to compute the groups in
(\ref{eq:FrNzerosection}) backwards, that is starting from $k=4$.

\begin{proposition}\label{prop:k=4}
${\rm H}^{i}({\bf G}/{\bf B},\wedge ^{4}{\mathcal
  T}\otimes {\sf F}_{\ast}{\sf S}^{\bullet}{\mathcal T}) = {\rm
  H}^{i}({\bf G}/{\bf B}, \omega _{{\bf G}/{\bf B}}^{-1} \otimes {\sf
  F}_{\ast}{\sf S}^{\bullet}{\mathcal T}) = 0$ for $i>0$.
\end{proposition}

\begin{proof}
Follows from Theorem \ref{th:KLTvantheorem}, the line bundle $\omega _{{\bf
    G}/{\bf B}}^{-1}$ being ample.
\end{proof}

\begin{proposition}\label{prop:k=4}
${\rm H}^{i}({\bf G}/{\bf B},\wedge ^{3}{\mathcal
  T}\otimes {\sf F}_{\ast}{\sf S}^{\bullet}{\mathcal T}) = 0$ for $i>1$.
\end{proposition}

\begin{proof}
There is an isomorphism $\wedge ^3{\mathcal T} = \Omega _{{\bf G}/{\bf
    B}}^{1}\otimes \omega _{{\bf G}/{\bf B}}^{-1}$. Consider the short
exact sequence dual to (\ref{eq:tangseqforproj_pi}) tensored with
$\omega _{{\bf G}/{\bf B}}^{-1}$:
\begin{equation}
0\rightarrow \pi ^{\ast}\Omega _{\Pp ^3}^{1}\otimes \omega _{{\bf G}/{\bf
    B}}^{-1}\rightarrow \wedge ^3{\mathcal T}\rightarrow \pi
    ^{\ast}\Oo _{\Pp ^3}(4)\rightarrow 0.
\end{equation}

Take the dual to the Euler sequence:
\begin{equation}\label{eq:Eulerpullback}
0\rightarrow \pi ^{\ast}\Omega _{\Pp ^3}^1\rightarrow {\sf
  V}^{\ast}\otimes \pi ^{\ast}\Oo _{\Pp ^3}(-1)\rightarrow \Oo _{{\bf G}/{\bf B}}\rightarrow 0.
\end{equation}

Tensoring the latter sequence with $\omega _{{\bf G}/{\bf B}}^{-1}\otimes {\sf F}_{\ast}{\sf S}^{\bullet}{\mathcal T}$ and
using Theorem \ref{th:KLTvantheorem}, we get:
\begin{equation}
{\rm H}^{i}({\bf G}/{\bf B},\wedge ^3{\mathcal T}\otimes {\sf F}_{\ast}{\sf S}^{\bullet}{\mathcal T}) = 0 \
  \mbox{for} \ i>1.
\end{equation}
\end{proof}

\begin{proposition}
${\rm H}^{i}({\bf G}/{\bf B},\wedge ^{2}{\mathcal
  T}\otimes {\sf F}_{\ast}{\sf S}^{\bullet}{\mathcal T}) = 0$ for $i>2$.
\end{proposition}

\begin{proof} From the sequence (\ref{eq:extsqoftangentseq})
  we see that it is sufficient to show  ${\rm H}^{i}({\bf
  G}/{\bf B},\pi ^{\ast}{\mathcal T}_{\Pp ^3}\otimes \Oo _{\pi}(2)\otimes {\sf F}_{\ast}{\sf S}^{\bullet}{\mathcal
  T}) = 0$ for $i>2$. Indeed, take the exterior square of the Euler
  sequence:
\begin{equation}
0\rightarrow \Oo _{{\bf G}/{\bf B}}\rightarrow {\V}\otimes \pi
^{\ast}\Oo _{\Pp ^3}(1)\rightarrow \wedge ^2{\V}\otimes \pi ^{\ast}\Oo
_{\Pp ^3}(2)\rightarrow \wedge ^2\pi ^{\ast}{\mathcal T}_{\Pp ^3}\rightarrow 0.
\end{equation}

Tensoring this sequence with ${\sf F}_{\ast}{\sf S}^{\bullet}{\mathcal
  T}$ and using Theorem \ref{th:KLTvantheorem}, we get ${\rm
  H}^{i}({\bf G}/{\bf B},\wedge ^2\pi ^{\ast}{\mathcal T}_{\Pp
  ^3}\otimes {\sf F}_{\ast}{\sf S}^{\bullet}{\mathcal T}) = 0$ for
  $i>0$. Therefore, the statement will hold if we show ${\rm H}^{i}({\bf G}/{\bf B},\Oo
  _{\pi}(2)\otimes \pi ^{\ast}\Oo _{\Pp ^3}(1)\otimes {\sf F}_{\ast}{\sf
  S}^{\bullet}{\mathcal T}) = 0$ for $i\geq 3$. This vanishing is
  equivalent to ${\rm H}^{i}({\bf G}/{\bf B},\Oo _{\pi}(2p)\otimes \pi
  ^{\ast}\Oo _{\Pp ^3}(p)) = 0$ for $i\geq 3$ and ${\rm H}^{4}({\bf G}/{\bf B},{\sf
  F}_{\ast}\Omega _{{\bf G}/{\bf B}}^1\otimes \Oo _{\pi}(2)\otimes \pi
  ^{\ast}\Oo _{\Pp ^3}(1)) = 0$. The vanishing of the former group for $i\geq 3$
  follows from Lemma \ref{lem:linebundlescohvan} below. The latter
  group will vanish, as the dual to the Euler sequence shows, if ${\rm H}^{4}({\bf G}/{\bf B},\Oo _{\pi}(2p)\otimes \pi
  ^{\ast}\Oo _{\Pp ^3}(p-1)) = 0$ and ${\rm H}^{4}({\bf G}/{\bf B},\Oo _{\pi}(2p-2)\otimes \pi
  ^{\ast}\Oo _{\Pp ^3}(p)) = 0$. These are the top cohomology groups of line
  bundles on ${\bf G}/{\bf B}$. By Serre duality one has:
\begin{eqnarray}
& {\rm H}^{4}({\bf G}/{\bf B},\Oo _{\pi}(2p)\otimes \pi ^{\ast}\Oo
_{\Pp ^3}(p-1))={\rm H}^{0}({\bf G}/{\bf B},\Oo _{\pi}(-2p-2)\otimes \pi
^{\ast}\Oo _{\Pp ^3}(-p-3))^{\ast} \\
& {\rm H}^{4}({\bf G}/{\bf B},\Oo _{\pi}(2p-2)\otimes \pi
  ^{\ast}\Oo _{\Pp ^3}(p))={\rm H}^{0}({\bf G}/{\bf B},\Oo
  _{\pi}(-2p)\otimes \pi ^{\ast}\Oo _{\Pp ^3}(-p-4))^{\ast} &\nonumber 
\end{eqnarray}

The groups in the right hand side are the groups of global sections of
line bundles on ${\bf G}/{\bf B}$ that are both non-effective. Note
that the same argument gives ${\rm H}^{4}({\bf G}/{\bf B},\Oo _{\pi}(2p)\otimes \pi
  ^{\ast}\Oo _{\Pp ^3}(p)) = 0$. By the
Kempf vanishing these groups are zero, hence the statement.

\end{proof}

\begin{proposition}\label{prop:vanishforTi>1}
${\rm H}^{i}({\bf G}/{\bf B},{\mathcal T}\otimes {\sf F}_{\ast}{\sf S}^{\bullet}{\mathcal T}) = 0$ for $i>1$.
\end{proposition}
\begin{proof}
It follows from the Euler sequence and Theorem \ref{th:KLTvantheorem} that ${\rm H}^{i}({\bf G}/{\bf
  B},\pi ^{\ast}{\mathcal T}_{\Pp ^3}\otimes {\sf F}_{\ast}{\sf
  S}^{\bullet}{\mathcal T}) = 0$ for $i>0$. Therefore, it is sufficient to show
that 
\begin{equation}
{\rm H}^{i}({\bf G}/{\bf B}, \Oo _{\pi}(2)\otimes {\sf F}_{\ast}{\sf
  S}^{\bullet}{\mathcal T}) = 0 
\end{equation}

for $i\geq2$. By Lemma \ref{lem:KLTlemma} this will follow from the vanishing of cohomology groups:
\begin{equation}
{\rm H}^{i}({\bf G}/{\bf B},\Oo _{\pi}(2)\otimes {\sf F}_{\ast}\Omega
_{{\bf G}/{\bf B}}^j) = 0
\end{equation}

for $i>j+1$ and $j=0,1,2$. Let us treat each case separately. \par 
\noindent ($\bf i$): Let $j=0$. We need to prove
that ${\rm H}^{i}({\bf G}/{\bf B},\Oo _{\pi}(2p)) = 0$ for $i\geq
2$. This follows from Lemma \ref{lem:linebundlescohvan} below.\par 
\noindent ($\bf ii$): Let $j=1$. In this case we need to show that ${\rm H}^{i}({\bf
  G}/{\bf B},\Oo _{\pi}(2)\otimes {\sf F}_{\ast}\Omega _{{\bf G}/{\bf B}}^1) = 0$ for $i\geq 3$.
Take the dual to (\ref{eq:tangseqforproj_pi}):
\begin{equation}\label{eq:dualtotangseqforproj_pi}
0\rightarrow \pi ^{\ast}\Omega _{\Pp ^3}^{1}\rightarrow \Omega _{{\bf G}/{\bf
  B}}^1\rightarrow \Oo _{\pi}(-2)\rightarrow 0.
\end{equation}

Applying the functor ${\sf F}_{\ast}$ to this sequence and tensoring it
with $\Oo _{\pi}(2)$, we get:
\begin{equation}
0\rightarrow {\sf F}_{\ast}\pi ^{\ast}\Omega _{\Pp ^3}^{1}\otimes
  \Oo _{\pi}(2)\rightarrow {\sf F}_{\ast}\Omega _{{\bf G}/{\bf
  B}}^1\otimes \Oo _{\pi}(2)\rightarrow {\sf F}_{\ast}\Oo
  _{\pi}(-2)\otimes \Oo _{\pi}(2)\rightarrow 0.
\end{equation}

Apply \ ${\sf F}_{\ast}$ to the sequence (\ref{eq:Eulerpullback}) \ and then
tensor it with $\Oo _{\pi}(2)$. Keeping in mind that ${\rm H}^{i}({\bf
  G}/{\bf B},\Oo _{\pi}(2p))=0$ for $i\geq 2$ (by Lemma \ref{lem:linebundlescohvan}) , we see that if ${\rm H}^{i}({\bf
  G}/{\bf B},\Oo _{\pi}(2p)\otimes \pi ^{\ast}\Oo (-1)) = 0$ for
$i\geq 3$, then ${\rm H}^{i}({\bf G}/{\bf B},{\sf F}_{\ast}\pi
^{\ast}\Omega _{\Pp ^3}^{1}\otimes \Oo _{\pi}(2)) = 0$ for $i\geq
3$. Together with the vanishing of ${\rm H}^{i}({\bf G}/{\bf B},{\sf F}_{\ast}\Oo
  _{\pi}(-2)\otimes \Oo _{\pi}(2))={\rm H}^{i}({\bf
  G}/{\bf B},\Oo _{\pi}(2p-2))$ for $i\geq 3$ this will imply ${\rm H}^{i}({\bf
  G}/{\bf B},\Oo _{\pi}(2)\otimes {\sf F}_{\ast}\Omega _{{\bf G}/{\bf
    B}}^1) = 0$ for $i\geq 3$.  Lemma \ref{lem:linebundlescohvan} states that ${\rm H}^{i}({\bf G}/{\bf B},\Oo
_{\pi}(2p)\otimes \pi ^{\ast}\Oo (-1)) = 0$ for $i\geq 3$ and ${\rm H}^{i}({\bf
  G}/{\bf B},\Oo _{\pi}(2p-2)) = 0$ for $i\geq 3$.\par 

\noindent ($\bf iii$): Let $j=2$. We \ only \ need \ to \ show \ the \
vanishing of the
top cohomology group ${\rm H}^4({\bf G}/{\bf B},{\sf F}_{\ast}\Omega
_{{\bf G}/{\bf B}}^2\otimes \Oo _{\pi}(2)) = 0$. There is a short
exact sequence:
\begin{equation}
0\rightarrow \pi ^{\ast}\Omega _{\Pp ^3}^2\rightarrow \Omega _{{\bf
    G}/{\bf B}}^2\rightarrow \pi ^{\ast}\Omega _{\Pp ^3}^1\otimes \Oo
    _{\pi}(-2)\rightarrow 0.
\end{equation}

Applying the functor ${\sf F}_{\ast}$ to this sequence and tensoring it
with $\Oo _{\pi}(2)$, we get:
\begin{equation}
0\rightarrow {\sf F}_{\ast}\pi ^{\ast}\Omega _{\Pp ^3}^2\otimes \Oo
    _{\pi}(2)\rightarrow {\sf F}_{\ast}\Omega _{{\bf
    G}/{\bf B}}^2\otimes \Oo _{\pi}(2)\rightarrow {\sf F}_{\ast}(\pi ^{\ast}\Omega _{\Pp ^3}^1\otimes \Oo
    _{\pi}(-2))\otimes \Oo _{\pi}(2)\rightarrow 0.
\end{equation}

Consider the Koszul resolution:
\begin{equation}
0\rightarrow \pi ^{\ast}\Omega _{\Pp ^3}^2\rightarrow \wedge ^2{\sf
  V}^{\ast}\otimes \pi ^{\ast}\Oo (-2)\rightarrow {\sf
  V}^{\ast}\otimes \pi ^{\ast}\Oo (-1)\rightarrow \Oo _{{\bf G}/{\bf
  B}}\rightarrow 0.
\end{equation}

Applying to this resolution the functor ${\sf F}_{\ast}$ and tensoring it with
$\Oo _{\pi}(2)$, we see that the vanishing of ${\rm H}^{4}({\bf
  G}/{\bf B},\Oo _{\pi}(2p)\otimes \pi ^{\ast}\Oo (-2))$, and that of
the cohomology groups of line bundles in ($\bf i$) and ($\bf ii$),
will imply ${\rm H}^{4}({\bf G}/{\bf B},{\sf F}_{\ast}\pi ^{\ast}\Omega
_{\Pp ^3}^2\otimes \Oo _{\pi}(2))=0$. Finally, the top cohomology
group of the bundle ${\sf F}_{\ast}(\pi ^{\ast}\Omega _{\Pp ^3}^1\otimes \Oo
_{\pi}(-2))\otimes \Oo _{\pi}(2)$ will vanish provided that ${\rm H}^{4}({\bf G}/{\bf B},\Oo
_{\pi}(2p-2)\otimes \pi ^{\ast}\Oo (-1)) = 0$ (tensor the Euler sequence with $\Oo
_{\pi}(-2)$, apply the functor ${\sf F}_{\ast}$, and finally
tensor the obtained sequence with $\Oo _{\pi}(2)$). Arguing as above, we immediately get 
the vanishing of top cohomology groups ${\rm H}^{4}({\bf
  G}/{\bf B},\Oo _{\pi}(2p)\otimes \pi ^{\ast}\Oo (-2))$ and ${\rm H}^{4}({\bf G}/{\bf B},\Oo
_{\pi}(2p-2)\otimes \pi ^{\ast}\Oo (-1))$.
\end{proof} 

Finally, Lemma \ref{lem:linebundlescohvan} completes the proof of Theorem \ref{th:B_2}.

\end{proof}

\begin{lemma}\label{lem:linebundlescohvan}
The following cohomology groups are zero:
\begin{itemize}

\item [(i)] ${\rm H}^{i}({\bf G}/{\bf B},\Oo _{\pi}(2p))=0$ for $i\geq 2$,

\item [(ii)] ${\rm H}^{i}({\bf G}/{\bf B},\Oo _{\pi}(2p-2))=0$ for $i\geq 3$,

\item [(iii)] ${\rm H}^{i}({\bf G}/{\bf B},\Oo _{\pi}(2p)\otimes \pi
  ^{\ast}\Oo (-1))=0$  for $i\geq 3$,

\item [(iv)] ${\rm H}^{i}({\bf G}/{\bf B},\Oo _{\pi}(2p)\otimes \pi
  ^{\ast}\Oo (p)) = 0$ for $i\geq 3$,

\end{itemize}
\end{lemma}

\begin{proof}
We first note that cohomology of line bundles on the flag varieties of
groups of rank two have been thoroughly studied (see, for instance,
\cite{An}). One can prove a large part of Lemma
\ref{lem:linebundlescohvan} using Andersen's criterion for
(non)--vanishing of the first cohomology group of a line bundle
(Theorem \ref{th:Andth}). However, there are minor errors about vanishing behaviour of
${\rm H}^2$ in \cite{An}; to make the exposition transparent and
self-contained we explicitly show all the vanishings listed above.\\

\noindent $(i)$ Consider the short exact sequence on ${\bf G}/{\bf B}$:
\begin{equation}\label{eq:fundseqonX}
0\rightarrow \pi ^{\ast}\Oo _{\Pp ^3}(-1)\rightarrow q^{\ast}\Uu _2\rightarrow \Oo
_{\pi}(-1)\rightarrow 0.
\end{equation}

Taking the determinants, one gets $q^{\ast}\Oo _{{\sf Q}_3}(1)=\pi
^{\ast}\Oo _{\Pp ^3}(1)\otimes \Oo _{\pi}(1)$. Dualizing (\ref{eq:fundseqonX}) and applying ${\sf F}^{\ast}$ to it, we obtain:
\begin{equation}\label{eq:decompforF^*U^*}
0\rightarrow \Oo _{\pi}(p)\rightarrow q^{\ast}{\sf F}^{\ast}\Uu _2
^{\ast}\rightarrow \pi ^{\ast}\Oo _{\Pp ^3}(p)\rightarrow 0.
\end{equation}

Finally, tensoring the above sequence with $\Oo _{\pi}(p)$, one has:
\begin{equation}\label{eq:seqforOo_pi(2p)}
0\rightarrow \Oo _{\pi}(2p)\rightarrow q^{\ast}{\sf F}^{\ast}\Uu _2
^{\ast}\otimes \Oo _{\pi}(p)\rightarrow q^{\ast}\Oo _{{\sf Q}_ 3}(p)\rightarrow 0.
\end{equation}

Consider the group ${\rm H}^{i}({\bf G}/{\bf B},q^{\ast}{\sf F}^{\ast}\Uu _2
^{\ast}\otimes \Oo _{\pi}(p))$. One has ${\rm R}^{\bullet}q_{\ast}\Oo
_{\pi}(p)={\sf S}^{p-2}\Uu _2^{\ast}(1)[-1]$. Indeed, $\Oo
_{\pi}(p)=\pi ^{\ast}\Oo _{\Pp ^3}(-p)\otimes q^{\ast}\Oo _{{\sf
    Q}_3}(p)$. By relative Serre duality we get ${\rm
  R}^{1}q_{\ast}\pi ^{\ast}\Oo _{\Pp ^3}(-p)={\sf S}^{p-2}\Uu _2(-1)$.
The above isomorphism follows from the projection formula and the isomorphism ${\sf S}^k\Uu
_2\otimes \Oo _{{\sf Q}_3}(k)={\sf S}^{k}\Uu _2^{\ast}$. Therefore,
\begin{equation}
{\rm H}^{i+1}({\bf G}/{\bf B},q^{\ast}{\sf F}^{\ast}\Uu _2
^{\ast}\otimes \Oo _{\pi}(p)) = {\rm H}^{i}({\sf Q}_3,{\sf F}^{\ast}\Uu _2
^{\ast}\otimes {\sf S}^{p-2}\Uu _2^{\ast}(1)).
\end{equation}

Recall the short exact sequence (cf. the sequence (\ref{eq:seqrank2})
from Section \ref{subsec:P1bundles}):
\begin{equation}
0\rightarrow {\sf F}^{\ast}\Uu _2^{\ast}\rightarrow {\sf S}^p\Uu
_2^{\ast}\rightarrow {\sf S}^{p-2}\Uu _2^{\ast}(1)\rightarrow 0.
\end{equation}

Tensoring it with ${\sf F}^{\ast}\Uu _2^{\ast}$, we obtain:
\begin{equation}
0\rightarrow {\sf F}^{\ast}(\Uu _2^{\ast}\otimes \Uu _2^{\ast})\rightarrow {\sf S}^{p}\Uu
_2^{\ast}\otimes {\sf F}^{\ast}\Uu _2^{\ast}\rightarrow {\sf S}^{p-2}\Uu
_2^{\ast}(1)\otimes {\sf F}^{\ast}\Uu _2^{\ast}\rightarrow 0.
\end{equation}

Consider the middle term of this sequence. Using Proposition \ref{prop:trivialmapS^2toD^2} (with $\Ee = \Uu _2$
and $n=1$) we get the short exact sequence:
\begin{equation}
0\rightarrow q^{\ast}\Oo _{{\sf Q}_3}(p)\rightarrow {\sf S}^{p}\Uu
_2^{\ast}\otimes {\sf F}^{\ast}\Uu _2^{\ast}\rightarrow {\sf
  S}^{2p}\Uu _2^{\ast}\rightarrow 0.
\end{equation}

Thus, ${\rm H}^{i}({\sf Q}_3,{\sf S}^{p}\Uu
_2^{\ast}\otimes {\sf F}^{\ast}\Uu _2^{\ast})=0$ for $i>0$. Let us
show that ${\rm H}^i({\sf Q}_3,{\sf F}^{\ast}(\Uu
_2^{\ast}\otimes \Uu _2^{\ast})) = 0$ for $i>1$. Note that ${\sf F}^{\ast}(\Uu
_2^{\ast}\otimes \Uu _2^{\ast}) = {\sf F}^{\ast}{\sf S}^{2}\Uu
_2^{\ast}\oplus {\sf F}^{\ast}\wedge ^2\Uu _2^{\ast}$ (the tensor
square of $\Uu _2^{\ast}$ splits into the direct sum, the
characteristic $p$ is odd). The latter
bundle is isomorphic to $\Oo _{{\sf Q}_3}(p)$, hence its higher
cohomology vanish. As for the former, note that ${\mathcal
  T}_{{\sf Q}_3}={\sf S}^{2}\Uu _2^{\ast}$ on ${\sf Q}_3$. Consider the
  adjunction sequence for the embedding $i: {\sf Q}_3\hookrightarrow \Pp ({\sf W})$:
\begin{equation}
0\rightarrow {\mathcal T}_{{\sf Q}_3}\rightarrow i^{\ast}{\mathcal
  T}_{\Pp ({\sf W})}\rightarrow \Oo _{{\sf Q}_3}(2)\rightarrow 0.
\end{equation}

Applying ${\sf F}^{\ast}$ to this sequence and using ${\rm H}^{i}({\sf
  Q}_3,{\sf F}^{\ast}i^{\ast}{\mathcal T}_{\Pp ({\sf W})})=0$ for $i>0$
  (use the Euler sequence on $\Pp ({\sf W})$, restrict it to ${\sf Q}_3$
  and then apply ${\sf F}^{\ast}$), we get ${\rm H}^{i}({\sf Q}_3,{\sf
  F}^{\ast} {\mathcal T}_{{\sf Q}_3})=0$ for $i>1$. This gives  ${\rm H}^{i}({\sf Q}_3,{\sf F}^{\ast}(\Uu
_2^{\ast}\otimes \Uu _2^{\ast})) = 0$ for $i>1$, and, therefore, ${\rm
  H}^{i}({\bf G}/{\bf B},\Oo _{\pi}(2p))=0$ for $i\geq 2$.\par

$(ii)$. The group ${\rm H}^{i}({\bf G}/{\bf B},\Oo _{\pi}(2p-2))$ for
$i\geq 3$. Similarly, tensoring the sequence
(\ref{eq:decompforF^*U^*}) with $\Oo _{\pi}(-2)$, we obtain:
\begin{equation}
0\rightarrow \Oo _{\pi}(2p-2)\rightarrow q^{\ast}{\sf F}^{\ast}\Uu _2
^{\ast}\otimes \Oo _{\pi}(p-2)\rightarrow q^{\ast}\Oo _{{\sf Q}_
  3}(p-2)\otimes \pi ^{\ast}\Oo _{\Pp ^3}(2)\rightarrow 0.
\end{equation}

As above, we get:
\begin{equation}
{\rm H}^{i+1}({\bf G}/{\bf B},q^{\ast}{\sf F}^{\ast}\Uu _2 ^{\ast}\otimes \Oo _{\pi}(p-2)) = {\rm H}^{i}({\sf Q}_3,{\sf F}^{\ast}\Uu _2
^{\ast}\otimes {\sf S}^{p-4}\Uu _2^{\ast}(1)),
\end{equation}

and, by the projection formula,
\begin{equation}
{\rm H}^{i}({\sf Q}_3,{\sf F}^{\ast}\Uu _2 ^{\ast}\otimes {\sf
  S}^{p-4}\Uu _2^{\ast}(1)) = {\rm H}^{i}({\bf G}/{\bf B},q^{\ast}{\sf F}^{\ast}\Uu _2
^{\ast}\otimes \pi ^{\ast}\Oo _{\Pp ^3}(p-2)\otimes q^{\ast}\Oo _{{\sf
  Q}_3}(1)).
\end{equation}

Tensoring the sequence (\ref{eq:decompforF^*U^*}) with $\pi ^{\ast}\Oo _{\Pp ^3}(p-2)\otimes q^{\ast}\Oo _{{\sf
  Q}_3}(1)$, we see that the bundle $q^{\ast}{\sf F}^{\ast}\Uu _2
^{\ast}\otimes \pi ^{\ast}\Oo _{\Pp ^3}(p-2)\otimes q^{\ast}\Oo _{{\sf
  Q}_3}(1)$ is an extension of two line
  bundles: $q^{\ast}\Oo _{{\sf Q}_3}(p+1)\otimes \pi ^{\ast}\Oo _{\Pp
  ^3}(-2)$ and $\pi ^{\ast}\Oo _{\Pp ^3}(2p-2)\otimes
  q^{\ast}\Oo _{{\sf Q}_3}(1)$. The latter is an effective line
  bundle, hence does not have higher cohomology by the Kempf
  vanishing. As for the former, one has:
\begin{equation}
{\rm H}^{i+1}({\bf G}/{\bf B},q^{\ast}\Oo _{{\sf Q}_3}(p+1)\otimes \pi
^{\ast}\Oo _{\Pp ^3}(-2)) = {\rm H}^i({\sf Q}_3,\Oo _{{\sf Q}_3}(p)),
\end{equation}

and the higher cohomology of the latter vanish as well. Thus, 
${\rm H}^{i}({\bf G}/{\bf B},\Oo _{\pi}(2p-2))$ for
$i\geq 3$.

$(iii)$ and $(iv)$. These groups can be treated in a similar way using
the sequence (\ref{eq:seqforOo_pi(2p)}) and the above arguments.
\end{proof}

\begin{remark}
{\rm There is a shorter proof of Proposition \ref{prop:vanishforTi>1}
  that can be, in fact, extended to another proof of Theorem \ref{th:B_2}. Note
  that ${\bf G}/{\bf B}$ is embedded into the product $\Pp ^3\times
  {\sf Q}_3$; denote $i$ this embedding. Consider the adjunction
  sequence:
\begin{equation}
0\rightarrow {\mathcal T}_{{\bf G}/{\bf B}}\rightarrow
  i^{\ast}{\mathcal T}_{\Pp ^3\times {\sf Q}_3}\rightarrow \pi
  ^{\ast}\Oo _{\Pp ^3}(1)\otimes q^{\ast}\Uu _2^{\ast}\rightarrow 0.
\end{equation}

Indeed, the variety ${\bf G}/{\bf B}$ can be represented as the zero
locus of a section of the bundle $\Oo _{\Pp ^3}(1)\boxtimes \Uu
_2^{\ast}$ on $\Pp ^3\times {\sf Q}_3$, hence the sequence. Applying
${\sf F}^{\ast}$ to it and tensoring with ${\sf S}^{\bullet}{\mathcal
  T}_{{\sf Q}_3}$, we get:
\begin{equation}
0\rightarrow {\sf F}^{\ast}{\mathcal T}_{{\bf G}/{\bf B}}\otimes {\sf S}^{\bullet}{\mathcal
  T}_{{\sf Q}_3}\rightarrow
  {\sf F}^{\ast}i^{\ast}{\mathcal T}_{\Pp ^3\times {\sf Q}_3}\otimes {\sf S}^{\bullet}{\mathcal
  T}_{{\sf Q}_3}\rightarrow \pi
  ^{\ast}\Oo _{\Pp ^3}(p)\otimes {\sf F}^{\ast}q^{\ast}\Uu _2^{\ast}\otimes {\sf S}^{\bullet}{\mathcal
  T}_{{\sf Q}_3}\rightarrow 0.
\end{equation}

Tensoring the sequence (\ref{eq:decompforF^*U^*}) with $\Oo _{\Pp
  ^3}(p)$, we get (cf. Proposition \ref{prop:trivialmapS^2toD^2}):
\begin{equation}\label{eq:decompforF^*U^*}
0\rightarrow q^{\ast}\Oo _{{\sf Q}_3}(p)\rightarrow \pi ^{\ast}\Oo
_{\Pp ^3}(p)\otimes q^{\ast}{\sf F}^{\ast}\Uu _2
^{\ast}\rightarrow \pi ^{\ast}\Oo _{\Pp ^3}(2p)\rightarrow 0.
\end{equation}

Tensoring it with ${\sf S}^{\bullet}{\mathcal T}_{{\sf Q}_3}$ and
using Theorem \ref{th:KLTvantheorem}, we obtain ${\rm H}^i({\sf Q}_3,\pi
  ^{\ast}\Oo _{\Pp ^3}(p)\otimes {\sf F}^{\ast}q^{\ast}\Uu _2^{\ast}\otimes {\sf S}^{\bullet}{\mathcal
  T}_{{\sf Q}_3})=0$ for $i>0$. Clearly, ${\rm H}^i({\sf Q}_3,{\sf F}^{\ast}i^{\ast}{\mathcal T}_{\Pp ^3\times {\sf Q}_3}\otimes {\sf S}^{\bullet}{\mathcal
  T}_{{\sf Q}_3})=0$ for $i>1$. Therefore, ${\rm H}^i({\sf Q}_3, {\sf F}^{\ast}{\mathcal T}_{{\bf G}/{\bf B}}\otimes {\sf S}^{\bullet}{\mathcal
  T}_{{\sf Q}_3})=0$ for $i>1$.

}
\end{remark}

\begin{corollary}
Let ${\bf G}/{\bf B}$ be the flag variety of type either ${\bf A}_2$
or ${\bf B}_2$. Assume that $p>3$ (respectively, $p>5$). Then ${\sf F}_{\ast}\Oo _{{\bf G}/{\bf B}}$ is a tilting bundle. 
\begin{proof}
First, by Lemma \ref{lem:Frobgeneratorforp>h}, the bundle ${\sf
  F}_{\ast}\Oo _{{\bf G}/{\bf B}}$ is a generator in $\Dd ^{b}({\bf G}/{\bf B})$. Secondly, Theorem
  \ref{th:SL_3} (respectively, Theorem \ref{th:B_2}) provide that the first condition of Definition \ref{def:tiltingdefinition} is
satisfied. Hence, for specified $p$, there is an equivalence of categories:
\begin{equation}
\Dd ^{b}({\bf G}/{\bf B})\simeq \Dd ^{b}({\rm End}({\sf F}_{\ast}\Oo _{{\bf G}/{\bf B}}) - \mbox{\rm mod})).
\end{equation}

\end{proof}
\end{corollary}


\section{Toric Fano varieties}\label{sec:toricFanovar}


Here we work out several examples of Fano toric varieties. 
\begin{lemma}
Let $X$ be the projective bundle $\Pp (\Oo \oplus \Oo (1))$ over $\Pp
^n$. Then $\Ext ^{i}({\sf F}_{\ast}\Oo _X,{\sf F}_{\ast}\Oo _X) = 0$
for $i>0$.
\end{lemma}

\begin{proof} 
We first consider the case $n>1$. Denote $\Ee = \Oo \oplus \Oo (1)$, and let $\pi \colon X\rightarrow
\Pp ^n$ be projection. Recall the short sequence (see Subsection \ref{subsec:P1bundles}):
\begin{equation}\label{eq:seqtoric}
0\rightarrow \pi ^{\ast}{\sf F}_{\ast}\Oo _{\Pp ^n} \rightarrow
{\sf F}_{\ast}\Oo _{X} \rightarrow \pi ^{\ast}({\sf F}_{\ast}({\tilde \Ee})\otimes
\mbox{det} \ (\Ee ^{\ast}))\otimes \Oo _{\pi}(-1)\rightarrow 0.
\end{equation}

Here ${\tilde \Ee} = {\sf S}^{p-2}\Ee\otimes \mbox{det} 
\Ee$. We first observe that the sequence (\ref{eq:seqtoric}) splits,
that is 
\begin{equation}
{\sf F}_{\ast}\Oo _X = \pi ^{\ast}{\sf F}_{\ast}\Oo _{\Pp ^n}\oplus 
\pi ^{\ast}({\sf F}_{\ast}{\tilde \Ee}\otimes \mbox{det}\Ee
^{\ast})\otimes \Oo _{\pi}(-1).
\end{equation}

Indeed, by adjunction one obtains an isomorphism 
\begin{equation}
\Ext ^{1}(\pi ^{\ast}({\sf F}_{\ast}{\tilde \Ee}\otimes \mbox{det}\Ee
^{\ast})\otimes \Oo _{\pi}(-1),\pi ^{\ast}{\sf F}_{\ast}\Oo _{\Pp ^n}) = 
\Ext ^{1}({\sf F}_{\ast}{\tilde \Ee}\otimes \mbox{det}\Ee
^{\ast},{\sf F}_{\ast}\Oo _{\Pp ^n}\otimes \pi _{\ast}\Oo _{\pi}(1)),
\end{equation}

the group in the right hand side being isomorphic to 
\begin{equation}\label{eq:Ext^1between2Frob}
\Ext ^{1}({\sf F}_{\ast}{\tilde \Ee},{\sf F}_{\ast}\Oo _{\Pp ^n}\otimes {\Ee}).
\end{equation}

Recall that ${\tilde \Ee} = {\sf S}^{p-2}\Ee\otimes \mbox{det}\Ee = \bigoplus _{k=1}^{k=p-1}\Oo (k)$.
It is well known that the Frobenius push-forward of any line
bundle on $\Pp ^n$ splits into direct sum of line bundles (e.g.,\cite{HKR}). More
generally, the Frobenius push-forward of a line bundle on a smooth
toric variety has the same property \cite{Bogv}. Hence, both terms
in the group (\ref{eq:Ext^1between2Frob}) are direct sums of line
bundles. However, a line bundle on $\Pp ^n$ does not have intermediate 
cohomology groups. Thus, the group
(\ref{eq:Ext^1between2Frob}) is zero and the bundle ${\sf F}_{\ast}\Oo _X$
splits.\\

We need to prove that $\Ext ^{i}(\pi
^{\ast}({\sf F}_{\ast}{\tilde \Ee}\otimes \mbox{det}\Ee 
^{\ast})\otimes \Oo _{\pi}(-1),{\sf F}_{\ast}\Oo _X) = 0$ for
$i>0$. This reduces to showing that 
\begin{equation}\label{eq:PO+O1vanishing}
\Ext ^{i}({\sf F}_{\ast}{\tilde \Ee},{\sf F}_{\ast}{\tilde \Ee})=0, \
\mbox{\rm and} \ \Ext ^{i}({\sf F}_{\ast}{\tilde \Ee},{\sf F}_{\ast}\Oo _{\Pp ^n}\otimes
{\Ee})=0
\end{equation}

for $i>0$.

\begin{proposition}\label{eq:FrobOlonP^n}
{\rm Let $1\leq k\leq p-1$. Then line bundles that occur in the decomposition of ${\sf F}_{\ast}\Oo
(k)$ are isomorphic to $\Oo (l)$ for $-n\leq l\leq 0$.}
\end{proposition}

\begin{proof}
Let ${\sf F}_{\ast}\Oo (k) = \oplus \ \Oo (a_i)$. Tensoring ${\sf
  F}_{\ast}\Oo (k)$ with $\Oo (-1)$, we obtain the bundle
${\sf F}_{\ast}\Oo (k)\otimes \Oo (-1) = {\sf F}_{\ast}\Oo (k-p)$, and
  the latter bundle has no global sections by the assumption. Hence,
  all $a_i\leq 0$. On the other hand, ${\sf F}_{\ast}\Oo (k)$ for such
  $k$ has no higher cohomology. This gives $a_i\geq -n$. 
\end{proof}

Lemma \ref{eq:FrobOlonP^n} implies that the groups in (\ref{eq:PO+O1vanishing})
are zero for $i>0$, hence the statement.\\

Now look at the case $n=1$. The toric variety -- the projective bundle $\Pp
(\Oo \oplus \Oo (1))$ over $\Pp ^1$ -- is a ruled surface ${\mathbb
  F}_1$ and is isomorphic to the blowup of a point on $\Pp
^2$. Blowups of $\Pp ^2$ are treated in the next section.
A straightforward check using Lemma \ref{cor:isomcor}, however, gives that the sequence
(\ref{eq:seqtoric}) splits for $n=1$ as well, and the rest of the
proof is the same as above.\par

In fact, the decomposition of the bundle ${\sf F}_{\ast}\Oo _X$ into a 
direct sum of line bundles allows to check when ${\sf F}_{\ast}\Oo _X$
generates the category $\Dd ^{b}(X)$. Let us treat the simplest case
of $\Pp ^1$. Recall the decomposition:
\begin{equation}
{\sf F}_{\ast}\Oo _{{\mathbb F}_1} = \pi ^{\ast}{\sf F}_{\ast}\Oo _{\Pp ^1}\oplus 
\pi ^{\ast}({\sf F}_{\ast}{\tilde \Ee}\otimes \mbox{det}\Ee
^{\ast})\otimes \Oo _{\pi}(-1).
\end{equation}

One has ${\sf F}_{\ast}\Oo _{\Pp ^1} = \Oo _{\Pp ^1}\oplus
\Oo _{\Pp ^1}(-1)^{\oplus p-1}$, and it can be easily verified (at
least for $p>2$) that a similar decomposition holds for ${\sf F}_{\ast}{\tilde
  \Ee}\otimes \mbox{det}\Ee ^{\ast}$:
\begin{equation}
{\sf F}_{\ast}{\tilde \Ee}\otimes \mbox{det}\Ee ^{\ast} =  \Oo _{\Pp
  ^1}^{\oplus a}\oplus
\Oo _{\Pp ^1}(-1)^{\oplus b},
\end{equation}

where $a$ and $b$ are non-zero multiplicities. It follows
  from Theorem \ref{th:Orlovth1} that the set of
  line bundles $\Oo _{{\mathbb F}_1},\pi ^{\ast}\Oo _{\Pp ^1}(-1),\Oo _
  {\pi}(-1), \pi ^{\ast}\Oo _{\Pp ^1}(-1)\otimes \Oo _{\pi}(-1)$
  generates the category $\Dd ^{b}({\mathbb F}_1)$, hence 
  for $p>2$ the bundle ${\sf F}_{\ast}\Oo _{{\mathbb F}_1}$ generates $\Dd
  ^{b}({\mathbb F}_1)$, that is ${\sf F}_{\ast}\Oo _{{\mathbb F}_1}$
  is a tilting bundle.
\end{proof}

Similarly, one checks the following:
\begin{lemma}
Let $X$ be the projective bundle $\Pp (\Oo \oplus \Oo (2))$ over $\Pp
^n$, and $n>2$. Then the bundle ${\sf F}_{\ast}\Oo _X$ is almost
exceptional (i.e. $\Ext ^{i}({\sf F}_{\ast}\Oo _X,{\sf F}_{\ast}\Oo _X) = 0$
for $i>0$).
\end{lemma}

\begin{proof}
The proof is completely analogous to that of the previous lemma. 
It turns out, however, that for $n=2$  there are non-vanishing
$\Ext$-groups in top degree. A direct sum decomposition for ${\sf
  F}_{\ast}\Oo _X$ as in the previous lemma holds anyway.
\end{proof}

\subsection{Toric Fano 3-folds.} Recall that according to the classification of smooth toric Fano
threefolds \cite{Ba}, there are 18 isomorphism classes of smooth
Fano toric threefolds. Among these are the following projective bundles: 
\begin{eqnarray}
& \Pp ^3, \ \Pp (\Oo _{\Pp ^2}\oplus \Oo _{\Pp ^2}(2)), \ \Pp (\Oo
  _{\Pp ^2}\oplus \Oo _{\Pp ^2}(1)), \ \Pp (\Oo _{\Pp ^1} \oplus \Oo
  _{\Pp ^1}\oplus \Oo _{\Pp ^1}(1)), \\
& \Pp (\Oo _{\Pp ^1\times \Pp ^1}\oplus \Oo _{\Pp ^1\times \Pp ^1}(1,1)), \
 \Pp (\Oo _{\Pp ^1\times \Pp ^1}\oplus \Oo _{\Pp ^1\times \Pp
  ^1}(1,-1)), \ \Pp (\Oo _{X_1}\oplus
  \Oo _{X_1}(l)), & \nonumber
\end{eqnarray}

\noindent The other varieties in the list are products of del Pezzo
  surfaces $X_k$ (the blowups of $\Pp ^2$ at $k$ points) and
  $\Pp ^1$:
\begin{equation}\label{eq:productsdelPezzotimesP^n}
\Pp ^2\times \Pp ^1, \Pp ^1\times \Pp ^1\times \Pp
  ^1,X_1\times \Pp ^1,X_2\times \Pp ^1,X_3\times \Pp ^1 ,
\end{equation}

and there are yet six varieties, of which four are
  isomorphic to del Pezzo fibrations over $\Pp ^1$.
The similar calculations as above give:\\

\noindent (i) Let $X$ be the projective bundle $\Pp (\Oo \oplus \Oo (1,1))$ over $\Pp
^1\times \Pp ^1$. Then $\Ext ^{i}({\sf F}_{\ast}\Oo _X,{\sf F}_{\ast}\Oo _X) = 0$
for $i>0$.\\

\noindent (ii) Let $X$ be the projective bundle $\Pp (\Oo \oplus \Oo (1,-1))$ over $\Pp
^1\times \Pp ^1$. Then $\Ext ^{1}({\sf F}_{\ast}\Oo _X,{\sf
  F}_{\ast}\Oo _X)\neq  0$.\\

In the next section we will see that $\Ext ^{i}({\sf F}_{\ast}\Oo
_{X_k},{\sf F}_{\ast}\Oo _{X_k}) = 0$ for a del Pezzo surface $X_k$. Thus,
the $\Ext$-groups vanish for all varieties in
(\ref{eq:productsdelPezzotimesP^n}).\\

Lemma \ref{lem:vanishviaomega} gives that for a Frobenius split variety $X$ the bundle ${\sf F}_{\ast}\Oo _X$ is almost exceptional
if the bundle ${\sf F}_{\ast}\Oo _X\otimes \omega _{X}^{-1}$ is
$\rm F$-ample. On the other hand, a line bundle is $\rm F$-ample if and only
if it is ample (\cite{Ar}, Lemma 2.4). By \cite{Bogv}, the bundle ${\sf F}_{\ast}\Oo _X\otimes
\omega _{X}^{-1}$ is the direct sum of line bundles (these line
bundles can explicitly be read off from the data defining the toric variety). 
\begin{proposition}\label{prop:easyprop}
Let $X$ be a smooth Fano toric variety. Conisder the decomposition 
$$
{\sf F}_{\ast}\Oo _X = \bigoplus_{i=0}^{i=N} \Ll _i
$$

\noindent where $\Ll _i$ are line bundles. If all $\Ll _i\otimes
\omega _X^{-1}$ are ample, then the bundle ${\sf F}_{\ast}\Oo _X$ is almost exceptional.
\end{proposition}

In the simplest example when $X=\Pp ^n$ we see that all the line
bundles in the decomposition are ample (use the argument from
Proposition \ref{eq:FrobOlonP^n}), hence ${\sf F}_{\ast}\Oo _{\Pp ^n}$
is almost exceptional. Corollary \ref{cor:Frobgeneratorforp>h}
then states that for $p>n+1$ the bundle ${\sf F}_{\ast}\Oo _{\Pp ^n}$
is tilting (cf. \cite{HKR}).


\section{Blowups of $\Pp ^2$}\label{sec:blowups}


\begin{theorem}\label{th:DelPezzotilting}
Let $X_k$ be a smooth surface that is obtained by blowing up of a set
of $k$ points on $\Pp ^{2}$ in general position, $k\geq 1$. Then for
$n\geq 1$ 
\begin{equation}
\Ext ^{i}({\sf F}^n_{\ast}\Oo _{X_k},{\sf F}^n_{\ast}\Oo _{X_k}) = 0
\end{equation}

for $i>0$.

\end{theorem}

\begin{proof}
We prove the theorem by induction, the base of induction being $X_0 =
 \Pp ^{2}$.  It is known that ${\sf F}_{\ast}\Oo _{\Pp ^2}$ is a generator in $\Dd ^{b}(\Pp ^2)$.
For $k\geq 1$ let $X_k = {\tilde \Pp ^2}_{x_1,\dots,x_{k}}$ be the blowup of $X_0$
at $k$ points in general position. Assume that 
\begin{equation}
\Ext ^{i}({\sf F}^n_{\ast}\Oo _{X_k},{\sf F}^n_{\ast}\Oo _{X_k}) = 0
\end{equation}

Consider the diagram:

\begin{figure}[H]
$$
\xymatrix @C5pc @R4pc {
 l_{k+1}\ar[r]^{i} \ar@<-0.1ex>[d]^{p} & X_{k+1} \ar@<-0.4ex>[d]^{\pi _k} \\
       x_{k+1}   \ar[r] & X_k
 }
$$
\end{figure}

Recall the short exact sequence:
\begin{equation}\label{eq:DelPezzoFrobseq}
0\rightarrow \pi _k^{\ast}{\sf F}^n_{\ast}\Oo _{X_k}\rightarrow {\sf F}^n_{\ast}\Oo
_{X_{k+1}}\rightarrow i_{\ast}\Oo
_{l_k}(-1)\otimes {\sf W}_n\rightarrow 0 .
\end{equation}

Applying the functor $\Hom ({\sf F}^n_{\ast}\Oo _{X_{k+1}}, ?)$ to it
, we get the long exact sequence:
\begin{eqnarray}
& 0\rightarrow \Hom ({\sf F}^n_{\ast}\Oo _{X_{k+1}},\pi _k^{\ast}{\sf F}^n_{\ast}\Oo
_{X_k})\rightarrow \Hom ({\sf F}^n_{\ast}\Oo _{X_{k+1}},{\sf F}^n_{\ast}\Oo
_{X_{k+1}})\rightarrow & \\
& \rightarrow \Hom ({\sf F}^n_{\ast}\Oo _{X_{k+1}},i_{\ast}\Oo
_{l}(-1)\otimes {\sf W}_n)\rightarrow \Ext ^{1}({\sf F}^n_{\ast}\Oo
_{X_{k+1}},\pi _k^{\ast}{\sf F}^n_{\ast}\Oo _{X_k})\rightarrow \dots
\quad .&\nonumber
\end{eqnarray}

Let us first consider the groups $\Ext ^{i}({\sf
  F}^n_{\ast}\Oo _{X_{k+1}},\pi _k^{\ast}{\sf F}^n_{\ast}\Oo
_{X_k})$. 

\begin{lemma}
$\Ext ^{i}({\sf F}^n_{\ast}\Oo _{X_{k+1}},\pi _k^{\ast}{\sf
  F}^n_{\ast}\Oo _{X_k}) = 0$ for $i>0$.
\end{lemma}

\begin{proof}
By adjunction we have:
\begin{eqnarray}\label{eq:1stgroupinDelPezzo}
& \Ext ^{i}({\sf F}^n_{\ast}\Oo _{X_{k+1}},\pi _k^{\ast}{\sf F}^n_{\ast}\Oo
_{X_k}) = \Ext ^{i}(\Oo _{X_{k+1}},{{\sf F}^{n}}^{\ast}\pi _k^{\ast}{\sf F}^n_{\ast}\Oo
_{X_k}\otimes \omega _{X_{k+1}}^{1-p^n}) = \\
& = {\rm H}^{i}(X_{k+1},\pi _k^{\ast}{{\sf F}^n}^{\ast}{\sf F}^n_{\ast}\Oo
_{X_k}\otimes \omega _{X_{k+1}}^{1-p^n}). &\nonumber
\end{eqnarray}

Recall that the canonical sheaves are related by the formula:
\begin{equation}\label{eq:cansheafforblowup}
\omega _{X_{k+1}} = \pi _k^{\ast}\omega _{X_k}\otimes \Oo _{X_{k+1}}(l_{k+1}).
\end{equation}

For any $m > 0$ there is a short exact sequence:
\begin{equation}\label{eq:canseqforblowup}
0\rightarrow \Oo _{X_{k+1}}(-ml_{k+1})\rightarrow \Oo _{X_{k+1}}\rightarrow
\Oo _{ml_{k+1}}\rightarrow 0.
\end{equation}

The sheaf $\Oo _{ml_{k+1}}$ has a filtration with associated graded
factors being ${\mathcal J}_{l_{k+1}}^{j}/{\mathcal J}_{l_{k+1}}^{j+1} = \Oo
_{l_{k+1}}(j)$ for $0\leq j < m$. Hence, ${\rm R}^{\bullet}{\pi _k}_{\ast}\Oo _{ml_{k+1}} = \Oo _{mx_{k+1}}$. 
Applying the direct image functor ${\pi _k}_{\ast}$ to the sequence (\ref{eq:canseqforblowup}), we get:
\begin{equation}
0\rightarrow {\pi _k}_{\ast}\Oo _{X_{k+1}}(-ml_{k+1})\rightarrow \Oo
_{X_k}\rightarrow \Oo _{mx_{k+1}}\rightarrow 0.
\end{equation}

Thus, ${\rm R}^{\bullet}{\pi _k}_{\ast}\Oo _{X_{k+1}}(-ml_{k+1}) = {\mathcal
  J}_{x_{k+1}}^{m}$ (cf. Theorem 5.7 in \cite{Har}). Finally:
\begin{equation}\label{eq:4147}
{\rm R}^{\bullet}{\pi _k}_{\ast}(\omega _{X_{k+1}}^{1-p^n}) = \omega _{X_k}^{1-p^n}\otimes
{\mathcal J}_{x_{k+1}}^{p^n-1}
\end{equation}

Using the projection formula, we obtain an isomorphism:
\begin{equation}
{\rm H}^{i}(X_{k+1},\pi _k^{\ast}{{\sf F}^n}^{\ast}{\sf F}^n_{\ast}\Oo
_{X_k}\otimes \omega _{X_{k+1}}^{1-p^n}) = 
{\rm H}^{i}(X_k,{{\sf F}^n}^{\ast}{\sf F}^n_{\ast}\Oo _{X_k}\otimes \omega _{X_k}^{1-p^n}\otimes
{\mathcal J}_{x_{k+1}}^{p^n-1}).
\end{equation}

Consider the short exact sequence 
\begin{equation}\label{eq:4155}
0\rightarrow {\mathcal J}_{x_{k+1}}^{p^n-1}\rightarrow \Oo _{X_k}\rightarrow \Oo
_{(p^n-1)x_{k+1}}\rightarrow 0,
\end{equation}

and its tensor product with the vector bundle
${{\sf F}^n}^{\ast}{\sf F}^n_{\ast}\Oo _{X_k}\otimes \omega _{X_k}^{1-p^n} =
{{\sf F}^n}^{!}{\sf F}^n_{\ast}\Oo _{X_k} : = \Ee _k$
\begin{equation}\label{eq:4150}
0\rightarrow {\mathcal J}_{x_{k+1}}^{p^n-1}\otimes \Ee _k\rightarrow \Ee
_k\rightarrow \Ee _k\otimes \Oo _{(p^n-1)x_{k+1}}\rightarrow 0.
\end{equation}

By the induction assumption we have 
\begin{equation}\label{eq:4151}
{\rm H}^{i}(X_k,\Ee _k) = \Ext ^{i}({\sf F}^n_{\ast}\Oo _{X_k},{\sf F}^n_{\ast}\Oo _{X_k})
= 0
\end{equation}

for $i>0$. 
\begin{proposition}\label{prop:46}
The map ${\rm H}^{0}(X_{k},\Ee _k)\rightarrow {\rm H}^{0}(X_k, \Ee _k\otimes \Oo
_{(p^n-1)x_{k+1}})$ is surjective.
\end{proposition}

\begin{proof}
We again proceed by induction. Let $k=1$. The bundle $\Ee _0 = {{\sf F}^n}^{\ast}{\sf F}^n_{\ast}\Oo _{\Pp ^{2}}\otimes \omega _{\Pp
  ^{2}}^{1-p}$ is isomorphic to $\Oo _{\Pp ^{2}}(3p^n-3)\oplus \Oo _{\Pp
  ^2}(2p^n-3)^{\oplus p_1}\oplus \Oo _{\Pp ^{2}}(p^n-3)^{\oplus p_2}$, where $p_1$ and $p_2$
are the multiplicities. In fact, these are computed to be
\begin{equation}
p_1 = \frac{(p^n-1)(p^n+4)}{2}, \quad p_2 = \frac{(p^n-1)(p^n-2)}{2} \quad .
\end{equation}

A dimension count gives that there is a surjection ${\rm H}^{0}(\Pp ^{2},\Ee _0)\twoheadrightarrow
{\rm H}^{0}(\Pp ^{2},\Ee _0\otimes \Oo _{(p^n-1)x_1})$. Indeed, the
dimension of the group ${\rm H}^{0}(\Pp ^{2},\Ee _0)$ is given by 
\begin{equation}\label{eq:dimencalculus}
\mbox{dim} \ {\rm H}^{0}(\Pp ^{2},\Ee _0) = \binom{3p^n-1}2 + \binom{2p^n-1}2
\cdot p_1 + \binom{p^n-1}2 \cdot p_2.
\end{equation}

On the other hand, the space ${\rm H}^{0}(\Pp ^{2},\Ee _0\otimes
\Oo _{(p^n-1)x_1})$ imposes $\frac{p^n(p^n-1)}{2}(1+p_1+p_2)$ conditions
and one sees that the dimension of this space is less than the
right-hand side in (\ref{eq:dimencalculus}). 
Moreover, the dimension count shows that there is even a surjection ${\rm H}^{0}(\Pp ^{2},\Ee _0)\twoheadrightarrow
{\rm H}^{0}(\Pp ^{2},\Ee _0\otimes \Oo _{p^nx_1})$. Fix now $k\geq 1$. Assume that for all $l\leq k$ we have a surjection 
${\rm H}^{0}(X_{l},\Ee _l)\rightarrow {\rm H}^{0}(X_l, \Ee _l\otimes \Oo
_{p^nx_{l+1}})$; in particular, this inductive assumption implies
a surjection ${\rm H}^{0}(X_{k},\Ee _k)\rightarrow {\rm H}^{0}(X_k, \Ee _k\otimes \Oo
_{(p^n-1)x_{k+1}})$, or, equivalently, ${\rm H}^{1}(X_k,{\mathcal J}_{x_{k+1}}^{p^n-1}\otimes
\Ee _k)=0$. 
Consider the sequence (\ref{eq:DelPezzoFrobseq}). Applying to it the functor
${{\sf F}^n}^{!}$ on $X_{k+1}$, we obtain:
\begin{equation}\label{eq:4158}
0\rightarrow {{\sf F}^n}^{!}\pi _k^{\ast}{\sf F}^n_{\ast}\Oo _{X_k}\rightarrow \Ee
_{k+1} \rightarrow {{\sf F}^n}^{!}(i_{\ast}\Oo
_{l_{k+1}}(-1)\otimes {\sf W}_n)\rightarrow 0.
\end{equation}

One has 
\begin{eqnarray}
& {{\sf F}^n}^{!}\pi _k^{\ast}{\sf F}^n_{\ast}\Oo _{X_k} =
{{\sf F}^n}^{\ast}\pi _k^{\ast}{\sf F}^n_{\ast}\Oo _{X_k}\otimes \omega
_{X_{k+1}}^{1-p^n} = \\
& = \pi _{k}^{\ast}({{\sf F}^n}^{\ast}{\sf F}^n_{\ast}\Oo _{X_k}\otimes
\omega _{X_k}^{1-p^n})\otimes \Oo _{X_{k+1}}(-(p^n-1)l_{k+1}) = 
\pi _{k}^{\ast}\Ee _k\otimes \Oo _{X_{k+1}}(-(p^n-1)l_{k+1}). & \nonumber
\end{eqnarray}

Denote the line bundle $\Oo _{X_{k+1}}(-(p^n-1)l_{k+1})$ by $\Ll _k$ and the sheaf 
${{\sf F}^n}^{!}(i_{\ast}\Oo _{l_{k+1}}(-1)\otimes {\sf W}_n)$ by ${\mathcal
  M}_k$. Tensor the sequence (\ref{eq:4158}) with the sheaf ${\mathcal
  J}_{x_{k+2}}^{p^n-1}$:
\begin{equation}\label{eq:4158tensorJ}
0\rightarrow {\mathcal Tor}^{1}({\mathcal M}_k,{\mathcal
  J}_{x_{k+2}}^{p^n-1})\rightarrow \pi _k^{\ast}\Ee _k\otimes \Ll
  _k\otimes {\mathcal J}_{x_{k+2}}^{p^n-1}
\rightarrow \Ee _{k+1} \otimes {\mathcal J}_{x_{k+2}}^{p^n-1}\rightarrow {\mathcal
  M}_k\otimes {\mathcal J}_{x_{k+2}}^{p^n-1} \rightarrow 0.
\end{equation}

We need to prove that ${\rm H}^{1}(X_{k+1},\Ee _{k+1}\otimes {\mathcal
  J}_{x_{k+2}}^{p^n-1})=0$. First observe that ${\mathcal
  M}_k\otimes {\mathcal J}_{x_{k+2}}^{p^n-1} = {\mathcal M}_k$ since the sheaf
  ${\mathcal J}_{x_{k+2}}$ is isomorphic to $\Oo _{X_{k+1}}$ in the
  formal neighborhood of the support of sheaf ${\mathcal M}_k$. Further,
for any coherent sheaf $\Ff$ on $X_{k+1}$ one has an exact sequence:
\begin{equation}\label{eq:H1forblowups}
0\rightarrow {\rm H}^{1}(X_k,{\rm R}^{0}{\pi _k}_{\ast}\Ff)\rightarrow
{\rm H}^{1}(X_{k+1},\Ff)\rightarrow {\rm H}^{0}(X_k,{\rm R}^{1}{\pi
  _k}_{\ast}\Ff)\rightarrow \dots
\end{equation} 

Applying to the short exact sequence 
\begin{equation}\label{eq:resolfornormalbundletoexceptcurve}
0\rightarrow \Oo _{X_{k+1}}\rightarrow \Oo _{X_{k+1}}(l_{k+1})\rightarrow
i_{\ast}\Oo _{l_{k+1}}(-1)\rightarrow 0
\end{equation}

the functor ${{\sf F}^n}^{!}$, we get:
\begin{equation}
0\rightarrow \omega _{X_{k+1}}^{1-p^n}\rightarrow \omega
_{X_{k+1}}^{1-p^n}\otimes \Oo _{X_{k+1}}(p^nl_{k+1})\rightarrow
{{\sf F}^n}^{!}(i_{\ast}\Oo _{l_{k+1}}(-1))\rightarrow 0.
\end{equation}

Take the direct image ${\pi _k}_{\ast}$:
\begin{equation}
0\rightarrow \omega _{X_k}^{1-p^n}\otimes {\mathcal
  J}_{x_{k+1}}^{p^n-1}\rightarrow \omega _{X_k}^{1-p^n}\rightarrow {\pi
  _k}_{\ast}{{\sf F}^n}^{!}(i_{\ast}\Oo _{l_{k+1}}(-1))\rightarrow 0.
\end{equation}

Therefore
\begin{equation}\label{eq:zerodirimageforF!exceptdiv}
{\rm R}^{0}{\pi _k}_{\ast}({{\sf F}^n}^{!}(i_{\ast}\Oo _{l_{k+1}}(-1))=\omega _{X_k}^{1-p^n}\otimes
  \Oo _{(p^n-1)x_{k+1}},
\end{equation}

and 
\begin{equation}\label{eq:firstdirimageforF!exceptdiv}
{\rm R}^{1}{\pi _k}_{\ast}({{\sf F}^n}^{!}(i_{\ast}\Oo _{l_{k+1}}(-1))=0.
\end{equation}

We see that ${\rm R}^{1}{\pi _k}_{\ast}({\mathcal
  M}_k\otimes {\mathcal J}_{x_{k+2}}^{p^n-1})=0$ and ${\rm R}^{0}{\pi
  _k}_{\ast}({\mathcal M}_k\otimes {\mathcal J}_{x_{k+2}}^{p^n-1})$ is a
skyscraper sheaf. By (\ref{eq:H1forblowups}) we get:
\begin{equation}
{\rm H}^{1}(X_{k+1},{\mathcal M}_k) = {\rm H}^{1}(X_{k+1},{\mathcal
 M}_k\otimes {\mathcal J}_{x_{k+2}}^{p^n-1})=0.
\end{equation}

Two following observations finish the proof: first, by induction assumption, one has
\begin{equation}
{\rm H}^{1}(X_{k+1},\pi _k^{\ast}\Ee _k\otimes \Ll _k\otimes {\mathcal
  J}_{x_{k+2}}^{p^n-1})={\rm H}^{1}(X_k,\Ee _k\otimes {\mathcal J}_{x_{k+1}\cup
{\pi _k}(x_{k+2})}^{p^n-1})=0.
\end{equation}

Secondly, the sheaf ${\mathcal Tor}^{1}({\mathcal M}_k, {\mathcal J}_{x_{k+2}}^{p^n-1})$
is a torsion sheaf supported on the exceptional divisor $l_{k+1}$; hence, 
${\rm H}^{2}(X_{k+1},{\mathcal Tor}^{1}({\mathcal M}_k, {\mathcal
  J}_{x_{k+2}}^{p^n-1}))=0$. Considering the spectral sequence
associated to the sequence (\ref{eq:4158tensorJ}) we get ${\rm H}^{1}(X_{k+1},\Ee _{k+1}\otimes
{\mathcal J}_{x_{k+2}}^{p^n-1})=0$, q.e.d.
\end{proof}

Taking into account (\ref{eq:4151}) and Proposition \ref{prop:46},
from the long exact cohomology sequence associated to (\ref{eq:4150}) 
we get:
\begin{equation}
{\rm H}^{i}(X_k,\Ee _k\otimes {\mathcal J}_{x_{k+1}}^{p^n-1}) = 0
\end{equation}

for $i>0$. Hence, the
left-hand side group in (\ref{eq:1stgroupinDelPezzo}) is zero for $i >
0$. 
\end{proof}

Now consider the groups $\Ext ^{i}({\sf F}^n_{\ast}\Oo _{X_{k+1}},i_{\ast}\Oo
_{l_{k+1}}(-1)\otimes {\sf W}_n) = {\rm H}^{i}(X_{k+1},{{\sf
  F}^n}^{!}(i_{\ast}\Oo _{l_{k+1}}(-1)\otimes {\sf W}_n))$. From (\ref{eq:zerodirimageforF!exceptdiv}) and
(\ref{eq:firstdirimageforF!exceptdiv}) one sees immediately that 
${\rm H}^{i}(X_{k+1},{{\sf F}^n}^{!}(i_{\ast}\Oo _{l_{k+1}}(-1)\otimes {\sf W_n})
= 0$ for $i > 0$. Finally, from (\ref{eq:1stgroupinDelPezzo}) we get $\Ext ^{i}({\sf F}^n_{\ast}\Oo
_{X_{k+1}},{\sf F}^n_{\ast}\Oo _{X_{k+1}}) = 0$ for $i > 0$, and the theorem
follows.

\end{proof}

\begin{corollary}
Let ${\rm D}_k = l_1 + \dots + \l_k$ be the exceptional divisor on
$X_k$. Then 
$$\Ext ^{i}({\sf F}^n_{\ast}\Oo _{X_k}({\rm D}_k),{\sf F}^n_{\ast}\Oo
_{X_k}({\rm D}_k)) = 0$$ 

for $i>0$.
\end{corollary}

\begin{proof}
Recall that ${{\pi _k}_{\ast}}\Oo _{X_{k+1}}({\rm D}_{k+1})=\Oo
_{X_{k}}({\rm D}_{k})$. Hence the
short exact sequence:
\begin{equation}\label{eq:DelPezzoFrobseqforO_Xl}
0\rightarrow \pi _k^{\ast}{\sf F}^n_{\ast}\Oo _{X_k}({\rm D}_k)\rightarrow {\sf F}^n_{\ast}\Oo
_{X_{k+1}}({\rm D}_{k+1})\rightarrow i_{\ast}\Oo
_{l_{k+1}}(-1)\otimes {\sf U}_n\rightarrow 0,
\end{equation}

where ${\sf U}_n$ is a vector space. The rest of the proof is
completely similar to that of Theorem \ref{th:DelPezzotilting}.
\end{proof}

\end{document}